\documentclass[amsfonts]{amsart}

\usepackage{amsmath}
\usepackage{amsfonts}
\usepackage{amsthm}
\usepackage{amssymb}
\theoremstyle{plain}
\newtheorem{theo}{Theorem}[section]
\newtheorem{lem}[theo]{Lemma}
\newtheorem{prop}[theo]{Proposition}
\newtheorem{claim}{Claim}
\theoremstyle{definition}
\newtheorem{definition}[theo]{Definition}
\theoremstyle{remark}
\newtheorem{rem}[theo]{Remark}
\newcommand{\pro}{\nabla u_1\cdot\nabla \overline{u}_0}
\newcommand{\dgr}{|\nabla (u_0-u_1)|^2}
\newcommand{\RR}{{\mathbb{R}}}
\newcommand{\CC}{{\mathbb C}}

\newcommand{\der}{\partial}
\newcommand{\ep}{\epsilon}

\newcommand{\om} {\Omega}

\newcommand{\dw}{{\delta W}}
\def\qed{\hfill$\square$\vspace{0.5cm}}    
\begin{document}

\title[Size estimates in the complex case]{Size estimates for the EIT problem with one measurement: the complex case}

\author[Beretta, Francini and Vessella.]{Elena~Beretta}
\address{ Dipartimento di Matematica ``G. Castelnuovo''
Universit\`a
di Roma ``La Sapienza'', Piazzale Aldo Moro 5, 00185 Roma, Italy}
\email{beretta@mat.uniroma1.it}
\author[]{Elisa~Francini}
\address{Dipartimento di Matematica ``U. Dini'', Universit\`a di Firenze,
Viale Morgagni 67A, 50134 Firenze, Italy}
\email{francini@math.unifi.it}
 \author[]{Sergio~Vessella}
 \address{Dipartimento di Matematica per le Decisioni, Universit\`a di Firenze,
 Via delle Pandette 9, 50127 Firenze, Italy}
\email{sergio.vessella@dmd.unifi.it}


\date{\today}

\keywords{Inverse boundary problems, size estimates, unique continuation}

\subjclass[2000]{Primary 35R30; Secondary 35R25, 35B60}

\date{}

\begin{abstract}
In this paper we estimate the size of a measurable inclusion in terms of power measurements for a single applied boundary current. This problem arises in medical imaging for the screening of organs (see \cite{G}).   For this kind of problem one has to deal mathematically with the complex conductivity (admittivity) equation. In this case we are able to establish, for certain classes of admittivities, lower and upper bounds of the measure of the inclusion in terms of the power measurements. A novelty of our result is that we are able to estimate also the volume of inclusions having part of its boundary in common with the reference body. Our analysis is based on the derivation of  energy bounds and of fine quantitative estimates of unique continuation for solutions to elliptic equations.\end{abstract}
\maketitle
\section{Introduction}\label{sec1}

In this paper we consider a mathematical problem arising in electrical impedance tomography (EIT), a nondestructive technique to determine electrical properties of a medium from measurements of voltages and currents on the boundary.

More precisely let $\Omega$ be the region occupied by a conducting medium and, at a fixed frequency $\omega$, consider  the complex valued admittivity function
\[
\gamma(x)=\sigma(x)+i\omega\epsilon(x),
\]
where $\sigma(x)$ represents the electrical conductivity at the point $x\in\Omega$ and $\epsilon(x)$ the electrical permittivity at a point $x\in\Omega$.

EIT leads to the inverse problem of the determination of the admittivity $\gamma$ from electrical measurements on $\partial\Omega$.
This technique has several applications in medical imaging, nondestructive testing of materials and geophysical prospection of the
underground.
We refer to the review paper by Borcea (\cite{Bo}) and to (\cite{CIN}) for a wide bibliography on relevant examples of applications. For a variational approach of the admittivity equation see \cite{CG}.

Relevant medical applications of EIT are for example breast cancer detection, (see for example \cite{CIN}) and screening of organs in transplantation surgery (\cite{G}).
In these particular situations one can assume $\gamma$ to be of the form
\[
\gamma=\gamma_0\chi_{\Omega\backslash D}+\gamma_1\chi_D,
\]
where $D\subset \Omega$ is a measurable subset of $\Omega$ and $\gamma_0\neq\gamma_1$.  $D$ represents either the cancerous tissue  or the degraded tissue which has different admittivity from the surrounding healthy one occupied by $\Omega\backslash D$. In particular in organs screening $D$ represents the region occupied by the degraded tissue imbedded in the healthy tissue and an important test to decide the quality of the organ is to give  an estimate of the size of $D$ in terms of boundary observations (\cite{G}).

Let us describe the mathematical problem:
let $\Omega\subset \RR^n$, $n\geq 2$ be a smooth, bounded domain and $D\subset \Omega$ be a measurable subset of $\Omega$. We denote by $\gamma_0$ and $\gamma_1$ the admittivities of  $\Omega\backslash D$ and $D$ respectively with
\[
\gamma_0=\sigma_0+i\epsilon_0 \text{ and }\gamma_1=\sigma_1+i\epsilon_1,
\]
(for simplicity we set $\omega=1$) and we assume that
\[
\sigma_0\geq c_0>0, \quad \sigma_1\geq c_0>0,
\]
this last condition corresponding to dissipation of energy, and let
\[
\gamma=\gamma_0\chi_{\Omega\backslash D}+\gamma_1\chi_D.
\]
Let $h\in H^{-1/2}(\der\om)$ be a complex valued boundary current flux and consider the so called \textit{background} potential $u_0\in H^1(\om)$  generated by the flux $h$, solution of
\begin{eqnarray*}
    \left\{\begin{array}{rcl}
    \mbox{div}(\gamma_0\nabla u_0)&=&0\mbox{ in }\om,\\
    \gamma_0\frac{\der u_0}{\der \nu}&=& h\mbox{ su }\der\om,
    \end{array}\right.
\end{eqnarray*}
and let $u_1\in H^1(\om)$ the \textit{perturbed} potential generated by the flux $h$ in the presence of the inclusion $D$, solution to
\begin{eqnarray*}
    \left\{\begin{array}{rcl}
    \mbox{div}(\gamma\nabla u_1)&=&0\mbox{ in }\om,\\
    \gamma\frac{\der u_1}{\der \nu}&=& h\mbox{ su }\der\om,
    \end{array}\right.
\end{eqnarray*}
with some common normalization condition.

Consider now
\[
W_1=\int_{\der\om}h\overline{u}_1,
\]
which represents the power required to maintain the current $h$ in the presence of the inclusion $D$ and analogously define
\[W_0=\int_{\der\om}h\overline{u}_0,
\]
 the power required to maintain the current $h$ in the unperturbed medium. Let
 \[
 \delta W=W_1-W_0
 \]
be the so called \textit{power gap}.

We will show that, if the admittivities $\gamma_0$ and $\gamma_1$ are constant or if $\gamma_0$ and $\gamma_1$ are variable scalar admittivities with $\gamma_0$ satisfying $\Im\gamma_0\equiv 0$ and some extra conditions, then the measure of $D$, $|D|$, can be estimated in terms of $|\delta W|$. For, we follow the approach introduced in \cite{AR}, \cite{ARS}  who derived estimates of $|D|$ in terms of the power gap for the real conductivity equation.\\
A different approach to derive size estimates for real conductivity inclusions has been introduced in \cite{CV} when $D$ is made of several connected components each of small size. Here  the authors use multiple boundary measurements of particular form to derive optimal
asymptotic estimates of $D$. Recently Kang et al. (see \cite{KKM})  obtained sharp  bounds of the size of two dimensional conductivity inclusions from  a pair of  boundary measurements using classical variational principles.\\

We want to point out that in the the screening of organs it seems to be  crucial to consider complex admittivities since electrical permittivity plays an important role in discriminating between degraded and normal tissue (\cite{G}).\\
To derive our main results, as mentioned above, we follow the approach of \cite{AR}, \cite{ARS}  making  use of the following basic tools  :
\begin{itemize}
\item Energy bounds.
\item Quantitative estimates of unique continuation.\end{itemize}
More precisely the first step is to find energy bounds i.e. lower and upper bounds for $\int_D|\nabla u_0|^2$ in terms of  $|\delta W|$ and the second one is to  find lower and upper bounds for $\int_D|\nabla u_0|^2$ in terms of  $|D|$ by using regularity and quantitative estimates of unique continuation of solutions to elliptic equations. Unfortunately, differently from the conductivity case, the first step in the complex case seems not to work for arbitrary  admittivities but only for constant ones or for certain variable admittivities variable scalar admittivities (see assumption {\bf (H3)} in Section \ref{sec2}).\\
On the other hand we would like  to emphasize  that,   in \cite{AR} and \cite{ARS}, the authors make the following technical assumption
\[
d(D,\der\om)\geq d_0>0.
\]
Clearly this hypothesis is rather restrictive in the medical application we have in mind since regions of the degraded tissue might appear also at the surface of the organ.
In this paper  we remove this assumption and prove size estimates also for inclusions having part of their boundary in common with $\der\om$. This is accomplished deriving  fine quantitative estimates of unique continuation (Lemma 4.4),  using reflection principles and suitable changes of variables.

The paper is divided as follows: in Section \ref{sec2} we state our main assumptions and our main results.
In Section \ref{sec3} we derive  energy bounds of the form
\[
K_1|\delta W|\leq \int_D|\nabla u_0|^2\leq K_2|\delta W|.
\]
 In Section \ref{sec4} we list some useful tools concerning quantitative estimates of unique continuation.\\
Section \ref{sec5} is devoted to the proof of our main results. In particular we derive lower and upper bounds of the measure of the inclusion in terms of the energy of the background potential in $D$.
 Finally in the appendix (Section \ref{sec7}) we give, for the reader's convenience, the proof of the Doubling Inequality stated in Section \ref{sec4}.

\section{Main results}\label{sec2}
\subsection{Notation and main assumptions}\label{subsec2.1}
For every $x\in\RR^n$ let us set $x=(x^\prime,x_n)$ where $x^\prime\in\RR^{n-1}$ for $n\geq 2$.

Let $x_0\in\RR^n$ and $r>0$. We denote by $B_r(x_0)$ and $B_r^\prime(x^\prime_0)$ the open ball in $\RR^n$ centered at $x_0$ of radius $r$ and the open ball in $\RR^{n-1}$ centered at $x^\prime$ of radius $r$, respectively. We denote by $Q_l(x_0)=\{x\in \RR^n : |x_j-x_{0j}|\leq l, j=1,\cdots, n\}$ the cube with center at $x_0$ and side $2l$.\\
\begin{definition}
\label{Definition 2.1}(${C}^{k,1}$ regularity)
Let $\om$  be  a bounded domain  in $\RR^n$. Given $k$, with $k=0,1$, we say that $\partial\Omega$ or $\om$ is of ${\it class }$ $ C^{k,1}$
with constants  $ r_0,  M_0$, if, for any $P\in \partial\Omega$, there exists a rigid transformation of coordinates under which we have $P=0$ and

\[
\om\cap (B'_{r_0}(0)\times (-M_0r_0,M_0r_0))=\{x\in B'_{r_0}(0)\times (-M_0r_0,M_0r_0) :\, x_n> \psi(x^\prime)\},
\]
where $\psi$ is a $C^{k,1}$ function on $B^\prime_{r_0}(0)$ such that
\[\psi(0)=0,
\]
\[
|\nabla\psi(0)|=0,\quad \hbox {when } k =1,\]

\[
\|\psi\|_{C^{k,1}(B^\prime_{r_0})}\leq M_0 r_0.
\]
\end{definition}
For $z,w\in\CC^n$ we denote by $z\cdot w=\sum_{j=1}^nz_jw_j$.

\begin{rem}
  \label{rem:2.1}
  We use the convention to normalize all norms in such a way that their
  terms are dimensionally homogeneous with their argument and coincide with the
  standard definition when the dimensional parameter equals one.
  For instance, the norm appearing above is meant as follows when $k=1$
\begin{equation*}
  \|\psi\|_{{C}^{1,1}(B^\prime_{r_0})} =
   \|\psi\|_{{L}^{\infty}(B^\prime_{r_0})}+ r_0
  \|\nabla\psi\|_{{L}^{\infty}(B^\prime_{r_0})}+
  r_0^{2}|\nabla\psi|_{1, B^\prime_{r_0}},
\end{equation*}
where
\begin{equation*}
|\nabla\psi|_{1, B^\prime_{r_0}}= \sup_ {\overset{\scriptstyle x,
y\in B^\prime_{r_0}}{\scriptstyle x\neq y}}
\frac{|\nabla\psi(x)-\nabla\psi(y)|} {|x-y|}.
\end{equation*}

Similarly, given a function $u:\Omega\mapsto \CC$,
\begin{equation*}
\|u\|_{L^2(\Omega)}=r_0^{-1}\left(\int_\Omega |u|^2\right)
^{\frac{1}{2}},
\end{equation*}
\begin{equation*}
\|u\|_{H^1(\Omega)}=r_0^{-1}\left(\int_\Omega |u|^2
+r_0^2\int_\Omega|\nabla u|^2 \right)^{\frac{1}{2}},
\end{equation*}
and so on for boundary and trace norms such as
$\|\cdot\|_{H^{\frac{1}{2}}(\partial\Omega)}$,
$\|\cdot\|_{H^{-\frac{1}{2}}(\partial\Omega)}$.

\end{rem}

We denote by $\Omega_r$, $r>0$, the following set
$$
\Omega_r=\{x\in \Omega : \textrm{ dist} (x,\partial\Omega)>r\}.
$$

Let us now set our main assumptions.
\vskip 3truemm
\noindent{\bf (H1)} {\it Assumptions on} $\om$.

Let $M_0$, $M_1$, $r_0$ be positive numbers such that $M_0\geq1$.

We assume that
\begin{enumerate}
\item $\om$  is a bounded domain  in $\RR^n$ with connected boundary;
\item $\der\om$ has $C^{0,1}$ regularity with constants $r_0$ and $M_0$;
\item $|\om|\leq M_1r_0^n$.
\end{enumerate}

\vskip 2truemm

\noindent{\bf (H2)} {\it Assumptions on} $D$.

$D$ is a Lebesgue measurable subset of $\overline{\om}$ and
\vskip 2truemm

{\bf (H2a)}
there exists a positive constant $d_0$ such that $\mbox{dist}(D,\der\om)\geq d_0$

\noindent or

{\bf (H2b)}
there exist $r_1\in(0,r_0]$ and $P\in\partial\Omega$ such that
\[
D\subset \overline{\om}\backslash B_{r_1}(P).
\]

\vskip 2truemm

\noindent{\bf (H3)} {\it Assumptions on the coefficients}.

Let $c_0\in (0,1], \mu_0$ and $L$ be positive numbers.  We assume the reference medium and the inclusion  have admittivities $\gamma_0=\sigma_0+i\epsilon_0$ and $\gamma_1=\sigma_1+i\epsilon_1$ satisfying
\[
\sigma_j\geq c_0, \quad |\gamma_j|\leq c_0^{-1}\text{ in }\om, \text{ for }j=0,1
\]
and, moreover we assume that

{\bf (H3i)}
\[\gamma_0\text{ and }\gamma_1\text{ are  constants}\]
and we denote by
\[\mu_0=|\gamma_0-\gamma_1|>0\]
or

{\bf (H3ii)}
\[\ep_0(x)\equiv 0\text{ in }\om,\text{ and }|\sigma_0(x)-\sigma_0(y)|\leq \frac{L}{r_0}|x-y|\text{ for }x,y\in\om\]
and
\[ |\ep_1(x)|\geq\mu_0\quad\text{ or }\quad \sigma_1(x)-\sigma_0(x)\geq\mu_0\text{ in }\om. \]

\vskip 2truemm

\noindent{\bf (H4)} {\it Assumptions on the boundary data}.

\vskip 0.2truemm

{\bf (H4a)}
 $h\in H^{-1/2}(\der\om)$ be a complex valued nontrivial current density at $\der\om$ satisfying
\[
\int_{\der\om}h=0.
\]
or

{\bf (H4b)}
 $h\in H^{-1/2}(\der\om)$ be a complex valued nontrivial current density at $\der\om$ satisfying
\[
\int_{\der\om}h=0,
\]
and such that
\[
\textrm{supp }h\subset \Gamma_0:=\der\om\cap\overline{B_{r_1/2}(P)},
\]
for the same $r_1$ and $P$ as in assumption  {\bf (H2b)}.

\vskip 5truemm

Let us denote by $F(h)$ the frequency of  $h$, that is
\begin{equation}\label{freq}
F(h)=\frac{\left\Vert h \right\Vert
_{H^{-1/2}(\partial \Omega )}}{
\left\Vert h \right\Vert _{H^{-1}(\partial \Omega )}}.
\end{equation}

Let
\[
\gamma=\gamma_0\chi_{\Omega\backslash D}+\gamma_1\chi_D
\]
and let us consider the unique solution $u_1\in H^1(\om)$ of the problem 
\begin{equation}\label{4}
    \left\{\begin{array}{rcl}
    \mbox{div}(\gamma\nabla u_1)&=&0\mbox{ in }\om,\\
    \gamma\frac{\der u_1}{\der \nu}&=& h\mbox{ su }\der\om,\\
    \int_{\partial\Omega}u_1&=&0.
    \end{array}\right.
\end{equation}
Analogously we define the background potential $u_0\in H^1(\om)$ generated by the same current flux $h$, the unique solution to the problem
\begin{equation}\label{5}
    \left\{\begin{array}{rcl}
    \mbox{div}(\gamma_0\nabla u_0)&=&0\mbox{ in }\om,\\
    \gamma_0\frac{\der u_0}{\der \nu}&=& h\mbox{ su }\der\om,\\
    \int_{\partial\Omega}u_0&=&0.
    \end{array}\right.
\end{equation}

\subsection{The main theorems}\label{subsec2.2}
We shall denote by $W_1$ and $W_0$ the power necessary to maintain the current $h$ when $D$ is present or absent, respectively:
\[
W_1=\int_{\der\om}h\overline{u}_1=\int_{\Omega}\gamma\nabla u_1\nabla\overline{u}_1,
\]
and
\[
W_0=\int_{\der\om}h\overline{u}_0=\int_{\Omega}\gamma_0\nabla u_0\nabla\overline{u}_0.
\]
Let $\delta W=W_1-W_0$ be the \emph{power gap}.

We first state our main result in the case of inclusions $D$ strictly contained in $\om$.
\begin{theo}\label{teo2.1}
Let  $\Omega$ satisfy {\bf (H1)} and let  $D$ be a measurable subset of $\Omega$ satisfying  {\bf (H2a)}.  Let $\gamma_0$ and $\gamma_1$ satisfy  assumptions {\bf (H3)} and let $h$ satisfy assumption {\bf (H4a)}.

Then,
\[
C_1\left|\frac{\delta W}{W_0}\right|\leq \frac{|D|}{|\om|}\leq C_2\left|\frac{\delta W}{W_0}\right|^{1/p}
\]
where $C_1$ depends on the a priori parameters $ c_0, \mu_0, M_0, M_1, \frac{d_0}{r_0}, L$, and the numbers $p>1$ and $C_2$ depend on the same a priori parameters and, in addition, on $F(h)$.
\end{theo}
We now state our main result in the case of inclusions that might have part of the boundary in common with $\der\om$.
\begin{theo}\label{teo2.2}
Let  $\Omega$ satisfy {\bf (H1)} with $\partial\Omega\in C^{1,1}$ with costants $r_0, M_0$and let  $D$ be a measurable subset of $\Omega$ satisfying  {\bf (H2b)}.  Let $\gamma_0$ and $\gamma_1$ satisfy assumptions {\bf (H3)} and let $h$ satisfy assumption {\bf (H4b)}.

Then,
\[
C_1\left|\frac{\delta W}{W_0}\right|\leq \frac{|D|}{|\om|}\leq C_2\left|\frac{\delta W}{W_0}\right|^{1/p},
\]
where $C_1$ depends on the a priori parameters $ c_0, \mu_0, M_0, M_1, \frac{r_1}{r_0}, L$, and the numbers $p>1$ and $C_2$ depend on the same a priori parameters and, in addition, on $F(h)$.
\end{theo}
\section{Energy bounds}\label{sec3}
\subsection{Energy identities}\label{subsec3.1}
In this section, following the idea first introduced in \cite{KSS}, we use energy identities in order to derive suitable energy bounds.

Let  $\tilde{\gamma}$ be a complex admittivity  and let us define the sesquilinear form
\[
    a_{\tilde\gamma}(u,v)=\int_\om\tilde{\gamma}\nabla u\cdot\nabla \overline{v}.
\]
If  $u_{\tilde\gamma}$   is a solution to
\[
    \left\{\begin{array}{rcl}
    \mbox{div}(\tilde{\gamma}\nabla u_{\tilde\gamma})&=&0\mbox{ in }\om,\\
   \tilde{ \gamma}\frac{\der u_{\tilde\gamma}}{\der \nu}&=& h\mbox{ on }\der\om,
    \end{array}\right.
\]
then
\begin{equation}\label{8}
    a_{\tilde{\gamma}}(u_{\tilde{\gamma}},v)=\int_{\der\om}h\overline{v},\quad \forall v\in H^1(\om).
\end{equation}
We observe that  in general $a_{\tilde\gamma}$ is not complex symmetric,
\begin{eqnarray*}
a_{\tilde{\gamma}}(u,v)-a_{\tilde{\gamma}}(v,u)&=&\int_\om\tilde{\gamma}\left(\nabla u\cdot\nabla\overline{v}-\nabla v\cdot\nabla\overline{u}\right)\nonumber\\
&=&2i\int_\om\tilde{\gamma}\Im(\nabla u\cdot\nabla\overline{v}).
\end{eqnarray*}

\begin{lem}\label{identita}
Let $\gamma_0$ and  $\gamma_1$ in $L^{\infty}(\Omega)$, let $\gamma=\gamma_0\chi_{\Omega\backslash D}+\gamma_1\chi_D$
and let  $u_1$ and $u_0$ the solutions of (\ref{4}) and (\ref{5}) respectively.
The following identities hold:
\begin{equation}
\tag{{\bf id1}}\int_\om\gamma\left|\nabla(u_1-u_0)\right|^2-\int_D(\gamma_1-\gamma_0)
|\nabla u_0|^2=\delta W
+2i\int_\om\gamma\Im(\nabla u_1\cdot\nabla \overline{u}_0),
\end{equation}
\begin{equation}
\tag{{\bf id2}}\int_\om\gamma_0\left|\nabla(u_1-u_0)\right|^2+\int_D(\gamma_1-\gamma_0)|\nabla u_1|^2
=-\delta W
-2i\int_\om\gamma_0\Im(\nabla u_1\cdot\nabla \overline{u}_0),
\end{equation}
\begin{equation}
\tag{{\bf id3}}
\int_D(\gamma_0-\gamma_1)\nabla u_1\cdot\nabla \overline{u}_0=\delta W+2i\int_\om\gamma_0\Im(\nabla u_1\cdot\nabla \overline{u}_0),
\end{equation}
\begin{equation}
\tag{{\bf id4}}
\int_D(\gamma_1-\gamma_0)\nabla u_0\cdot\nabla \overline{u}_1=-\delta W-2i\int_\om\gamma\Im(\nabla u_1\cdot\nabla \overline{u}_0).
\end{equation}
\end{lem}
\textit{Proof.}
Let us denote by $a_0(u,v):=a_{\gamma_0}(u,v)$ and by  $a_1(u,v):=a_{\gamma}(u,v)$.

\noindent From (\ref{8}) we have
\[
    a_0(u_0,v)=a_1(u_1,v)=\int_{\der\om}h\overline{v},\quad\forall v\in H^1(\om).
\]
Let us compute
\begin{eqnarray}\label{J1a}
  J_1 &:=& a_1(u_1-u_0,u_1-u_0)-[a_1(u_0,u_0)-a_0(u_0,u_0)]\nonumber \\
   &=&\int_{\der\om}h\overline{u}_1-\int_{\der\om}h\overline{u}_0+2i\int_\om\gamma\Im(\nabla u_1\cdot\nabla \overline{u}_0).
\end{eqnarray}
On the other hand
\begin{eqnarray}\label{J1b}
J_1&=&\int_\om\gamma\left|\nabla(u_1-u_0)\right|^2-\int_\om(\gamma-\gamma_0)
|\nabla u_0|^2\nonumber\\&=&\int_\om\gamma\left|\nabla(u_1-u_0)\right|^2-\int_D(\gamma_1-\gamma_0)
|\nabla u_0|^2,
\end{eqnarray}
and so, by (\ref{J1a}) and (\ref{J1b}) and the definition of $\delta W$, identity {\bf(id1)} follows.

Analogously we can compute
\begin{eqnarray*}
 J_2 &:=& a_0(u_0-u_1,u_0-u_1)-[a_0(u_1,u_1)-a_1(u_1,u_1)] \\
   &=& -\int_{\der\om}h(\overline{u}_1-\overline{u}_0)-2i\int_\om\gamma_0\Im(\nabla u_1\cdot\nabla \overline{u}_0).
\end{eqnarray*}
On the other hand
\[
J_2=\int_\om\gamma_0\left|\nabla(u_1-u_0)\right|^2+\int_\om(\gamma_1-\gamma_0)
|\nabla u_1|^2
\]
and, hence, {\bf(id2)} follows.

Finally let us compute
\[a_0(u_1,u_0)-a_1(u_1,u_0)=\int_D(\gamma_0-\gamma_1)\nabla u_1\cdot\nabla \overline{u}_0,\]
and, observe that
\begin{eqnarray*}
a_0(u_1,u_0)-a_1(u_1,u_0)&=&a_0(u_1,u_0)-a_1(u_1,u_0)+a_0(u_0,u_1)-a_0(u_0,u_1)\\
&=&\int_{\der\om}h(\overline{u}_1-\overline{u}_0)+2i\int_\om\gamma_0\Im(\nabla u_1\cdot\nabla \overline{u}_0),
\end{eqnarray*}
so that {\bf(id3)} follows.

By symmetry we can also write {\bf (id4)}
\qed
\begin{rem}
Note that  by combining {\bf(id1)} and {\bf(id4)}, we get an easy consequence of the definition of $u_0$ and $u_1$, that is
\begin{equation}\label{loc}
\int_\om\gamma\left|\nabla(u_1-u_0)\right|^2=\int_D(\gamma_0-\gamma_1)\nabla(\overline{u}_1-\overline{u}_0)\nabla u_0.
\end{equation}
\end{rem}

\subsection{The constant case}\label{subsec3.3}
\begin{prop}\label{prop3.1}
Assume $\gamma_0$ and $\gamma_1$ satisfy ${\bf (H3i)}$ and let $u_0$ and $u_1$ be the solution of (\ref{5}) and (\ref{4}), then
\[
\frac{c_0}{(c_0+|\gamma_1-\gamma_0|)|\gamma_1-\gamma_0|}\,|\delta W|\leq \int_D|\nabla u_0|^2\leq \left(\frac{1}{c_0}+\frac{2}{|\gamma_1-\gamma_0|}\right)|\delta W|.
\]
\end{prop}

\textit{Proof. }
Since $\gamma_0$ is constant and not zero we can write
\begin{eqnarray*}
\int_\om  \Im(\nabla u_1\cdot\nabla\overline{u}_0)&=&-\int_\om\Im\left(\nabla u_0\cdot\nabla\overline{u}_1-\nabla u_0\cdot\nabla\overline{u}_0\right)\nonumber\\
&=& -\int_\om\Im\left(\gamma_0\left(\nabla u_0\cdot\nabla \overline{u}_1-\nabla u_0\cdot\nabla \overline{u}_0\right)\frac{1}{\gamma_0}\right)\nonumber\\
&=&-\Im\left(\frac{1}{\gamma_0}\int_\om\gamma_0\left(\nabla u_0\cdot\nabla \overline{u}_1-\nabla u_0\cdot\nabla \overline{u}_0\right)\right)\nonumber\\
&=&-\Im\left(\frac{1}{\gamma_0}\int_{\der\om}h(\overline{u}_1-\overline{u}_0)\right)
=-\Im\left(\frac{\dw}{\gamma_0}\right)
\end{eqnarray*}
and, hence,
\[  \int_\om  \gamma_0\Im(\nabla u_1\cdot\nabla\overline{u}_0)=-\gamma_0\Im\left(\frac{\dw}{\gamma_0}\right).\]

Then, if we set
\begin{equation}\label{deltaV}
\delta V=\delta W-2i\gamma_0\Im\left(\frac{\delta W}{\gamma_0}\right)=\delta W+2i\int_\om\gamma_0\Im(\nabla u_1\cdot\nabla \overline{u}_0)
\end{equation}
we can write the identities of Lemma \ref{identita} in the following way
\begin{equation}
\tag{{\bf id1c}}\int_\om\gamma\left|\nabla(u_1-u_0)\right|^2-\int_D(\gamma_1-\gamma_0)
|\nabla u_0|^2=2i\int_D(\gamma_1-\gamma_0)\Im\left(\nabla u_1\cdot\nabla \overline{u}_0\right)
+\delta V,
\end{equation}
\begin{equation}
\tag{{\bf id2c}}\int_\om\gamma_0\left|\nabla(u_1-u_0)\right|^2+\int_D(\gamma_1-\gamma_0)|\nabla u_1|^2
=-\delta V,\end{equation}
\begin{equation}
\tag{{\bf id3c}}
\int_D(\gamma_0-\gamma_1)\nabla u_1\cdot\nabla \overline{u}_0=\delta V,
\end{equation}
\begin{equation}
\tag{{\bf id4c}}
\int_D(\gamma_1-\gamma_0)\nabla u_0\cdot\nabla \overline{u}_1=-2i\int_D(\gamma_1-\gamma_0)\Im\left(\nabla u_1\cdot\nabla \overline{u}_0\right)-\delta V.
\end{equation}

Let us write
\begin{eqnarray}\label{graduzero}
\int_D|\nabla u_0|^2&=&\int_D\dgr-\int_D|\nabla u_1|^2+2\int_D\Re(\pro)\nonumber\\
&\leq &\int_{\Omega}\dgr-\int_D|\nabla u_1|^2+2\int_D\Re(\pro).
\end{eqnarray}
By taking the real part of {\bf (id2c)} we get
\[
\int_\om\sigma_0\dgr+(\sigma_1-\sigma_0)\int_D|\nabla u_1|^2=-\Re(\delta V);
\]
by dividing by the positive constant $\sigma_0$ and using the fact that both $\sigma_0$ and $\sigma_1$ are positive we have
\begin{equation}\label{13*}
\int_\om\dgr-\int_D|\nabla u_1|^2\leq -\frac{\Re(\delta V)}{\sigma_0}.
\end{equation}
Now, let us divide  {\bf (id3c)} by the constant $\gamma_0-\gamma_1\neq 0$ and take the real part. We get
\[
\int_D\Re(\pro)=\Re\left(\frac{\delta V}{\gamma_0-\gamma_1}\right),
\]
which, together with  (\ref{13*}) and (\ref{graduzero}), gives
\[
\int_D|\nabla u_0|^2
\leq -\frac{\Re( \delta V)}{\sigma_0}+2\Re\left(\frac{\delta V}{\gamma_0-\gamma_1}\right).\]
This leads to the following upper bound
\[
\int_D|\nabla u_0|^2
\leq |\delta V|\left(\frac{1}{c_0} +\frac{2}{|\gamma_0-\gamma_1|}\right).
\]

\noindent To prove the lower bound observe that by (\ref{loc}) and since $\Re \gamma\geq c_0$ we have
\begin{equation}\label{localizz}
\left(\int_\om\dgr\right)^{1/2} \leq\frac{|\gamma_0-\gamma_1|}{c_0}\left(\int_D|\nabla u_0|^2\right)^{1/2}.
\end{equation}

Hence using identity {\bf(id3c)} we have
\noindent
\begin{eqnarray*}
\left|\delta V\right|&=&\left|\int_D(\gamma_0-\gamma_1)\pro\right|\\
&=&\left|(\gamma_0-\gamma_1)\left(\int_D\nabla(u_1-u_0)\cdot\nabla\bar u_0+\int_D|\nabla u_0|^2\right)\right|\\
&\leq&|\gamma_0-\gamma_1|\left(\left(\int_D|\nabla (u_1-u_0)|^2\right)^{1/2}\left(\int_D|\nabla u_0|\right)^{1/2}+\int_D|\nabla u_0|^2\right)\\
&\leq&|\gamma_0-\gamma_1|\left(\frac{|\gamma_0-\gamma_1|}{c_0}\int_D|\nabla u_0|^2+\int_D|\nabla u_0|^2\right),
\end{eqnarray*}
from which the lower bound
\[\int_D|\nabla u_0|^2\geq \frac{1}{|\gamma_0-\gamma_1|\left(\frac{|\gamma_0-\gamma_1|}{c_0}+1\right)}\left|\delta V\right|\]
follows.

Now, by using (\ref{deltaV}), we can see that
\[\delta V=\frac{\gamma_0^2}{|\gamma_0|^2}\overline{\delta W}\]
hence, in particular,
\[|\delta V|=|\delta W|\]
and the thesis follows.\qed

\subsection{The variable case}
\begin{prop}\label{prop3.2}
Assume $\gamma_0$ and $\gamma_1$ satisfy assumption ${\bf (H3ii)}$ and let $u_0$ the solution of (\ref{5}), then
\begin{equation}\label{ebv}
K_1|\delta W|\leq \int_D|\nabla u_0|^2\leq K_2|\delta W|,
\end{equation}
where
\[K_1=\frac{c_0^3}{2(2+c_0^2)}\quad\mbox{ and }\quad K_2=2\left(\frac{1}{\mu_0c_0^2}+\frac{1}{\mu_0}+\frac{1}{c_0}\right).\]
\end{prop}

\textit{Proof}
If assumption ${\bf (H3ii)}$ holds, then $\gamma_0=\sigma_0$ and $\ep_0=0$.
In this case, we have

\begin{eqnarray*}
\int_\om \sigma_0 \Im(\nabla u_1\cdot\nabla\overline{u}_0)&=&\int_\om\sigma_0\Im\left(\nabla u_1\cdot\nabla\overline{u}_0-\nabla u_0\cdot\nabla\overline{u}_0\right)\nonumber\\
&=& \Im\left(\int_\om\sigma_0\nabla \overline{u}_0\cdot\nabla u_1-\sigma_0\nabla \overline{u}_0\cdot\nabla u_0\right)\nonumber\\
&=&\Im\left(\int_{\der\om}\overline{h}u_1-\int_{\der\om}\overline{h}u_0\right)=\Im(\overline{\dw})=-\Im(\dw),
\end{eqnarray*}
and the energy identities become

\begin{equation}
\tag{\bf id1*}\int_\om \gamma|\nabla(u_0-u_1)|^2-
\int_D(\gamma_1-\gamma_0)|\nabla u_0|^2= \overline{\dw} +2i\int_D(\gamma_1-\gamma_0)\Im(\nabla u_1\cdot\nabla \overline{u}_0),
\end{equation}
\begin{equation}
\tag{\bf id2*}\int_\om\gamma_0\dgr +\int_D(\gamma_1-\gamma_0)|\nabla u_1|^2=-\overline{\dw},
\end{equation}
\begin{equation}\tag{\bf id3*}
\int_D(\gamma_0-\gamma_1)\pro =\overline{\dw}.
\end{equation}

By identity {\bf (id3*)} we have that
\begin{eqnarray}\label{e1}\nonumber
|\overline{\delta W} |&=&\left|\int_D(\gamma_0-\gamma_1)\nabla u_1\nabla\overline{u}_0\right|\\&=&
\left|\int_D(\gamma_0-\gamma_1)\nabla(u_1-u_0)\nabla\overline{u}_0+\int_D(\gamma_0-\gamma_1)|\nabla u_0|^2\right|\nonumber\\
&\leq&\sup_D|\gamma_0-\gamma_1|\left(\left(\int_D\dgr\right)^{1/2}\left(\int_D|\nabla u_0|^2\right)^{1/2}+\int_D|\nabla u_0|^2\right).
\end{eqnarray}
By (\ref{loc}), we have
\[
\left(\int_\om|\nabla(u_1-u_0)|^2\right)^{1/2}\leq\frac{\sup_{D}|\gamma_0-\gamma_1|}{c_0}
\left(\int_D|\nabla u_0|^2\right)^{1/2},
\]
and by combining this with (\ref{e1}) we get
\[|\delta W|\leq \sup_D|\gamma_0-\gamma_1|\left(\frac{\sup_D|\gamma_0-\gamma_1|}{c_0}+1\right)\int_D|\nabla u_0|^2\leq \frac{2}{c_0}\left(\frac{2}{c_0^2}+1\right)\int_D|\nabla u_0|^2,\]
and one side of estimate (\ref{ebv}) follows.

To derive the upper bound, let us first assume
\begin{equation}\label{e2}
 \sigma_1-\sigma_0\geq\mu_0.
\end{equation}
\noindent From the real part of identity
{\bf (id2*)} we get
\begin{equation}\label{e2.5}
\int_\om\sigma_0|\nabla(u_1-u_0)|^2+\int_D(\sigma_1-\sigma_0)|\nabla u_1|^2=
-\Re(\dw),
\end{equation}
hence, by assumption (\ref{e2}),
\[\int_\om |\nabla(u_1-u_0)|^2\leq-\frac{\Re (\dw)}{c_0},\]
\[
   \int_D |\nabla u_1|^2\leq-\frac{\Re(\dw)}{\mu_0},
\]
so that
\[\int_D|\nabla u_0|^2\leq 2\int_D|\nabla (u_0-u_1)|^2+2\int_D|\nabla u_1|^2\leq-2\left(\frac{1}{c_0}+\frac{1}{\mu_0}\right)\Re(\dw).\]

On the other hand, if $|\ep_1|\geq \mu_0$, then , from the imaginary part of {\bf (id2*)}, we get
\[\int_D |\ep_1||\nabla u_1|^2=|\Im (\dw)|,\]
and, hence,
\begin{equation}\label{e4}
    \int_D|\nabla u_1|^2\leq \frac{|\Im(\dw)|}{\mu_0}.
\end{equation}
\noindent From the real part of {\bf (id2*)} (see (\ref{e2.5})) and from (\ref{e4}) we get
\begin{eqnarray}\label{e5}
  \int_\om|\nabla(u_1-u_0)|^2 &\leq& \int_\om\sigma_0 c_0^{-1}|\nabla(u_1-u_0)|^2 \nonumber\\
  &=& c_0^{-1} \int_D(\sigma_0-\sigma_1)|\nabla u_1|^2-c_0^{-1}\Re(\dw)\nonumber\\
   &\leq& c_0^{-1} \sup_D|\sigma_0-\sigma_1|\int_D|\nabla u_1|^2-c_0^{-1}\Re(\dw)\nonumber\\
   &\leq&  \frac{1}{c_0\mu_0} \sup_D|\sigma_0-\sigma_1|\,|\Im(\dw)|-c_0^{-1}\Re(\dw)\nonumber\\
   &\leq&\frac{1}{c_0^2\mu_0}|\Im(\dw)|-\frac{1}{c_0}\Re(\dw).
\end{eqnarray}
By (\ref{e4}) and (\ref{e5}) we get the upper bound
\[\int_D|\nabla u_0|^2\leq 2\left(\frac{1}{\mu_0c_0^2}+\frac{1}{c_0}+\frac{1}{\mu_0}\right)|\delta W|.\]
\qed

\subsection{A one dimensional example}
We are not able to derive  energy bounds and hence also estimates on the size of $D$ for arbitrary variable admittivities. On the other hand we have seen in Proposition \ref{prop3.2}  that assumption ${\bf (H3ii)}$ leads to energy estimates. The lack of symmetry of condition ${\bf (H3ii)}$, that seems not natural, is in some sense optimal as the following example shows.

Let  $\Omega= (-1,1)$ and let  $D=[a,b]\subset (-1,1)$. Consider the background solution $u_0$
\begin{eqnarray*}
 \left\{\begin{array}{rcl}
(\gamma_0 u'_0)' &=&0\text{ in }(-1,1),\\
(\gamma_0 u'_0) (-1)&=&(\gamma_0 u'_0) (1)=K\in \mathbb C, \quad u_0(-1)+u_0(1)=0.
 \end{array}\right.
 \end{eqnarray*}
Integrating the equation and using the normalization conditions one gets that
 \[
 u_0(x)=F_0(x)+M,\mbox{ for } x\in (-1,1),
 \]
where
\[
F_0(x)=\int \frac{K}{\gamma_0(x)} dx,
\]
and $M-\frac{(F_0(1)+F_0(-1))}{2}$.

\noindent Considering  the perturbed solution $u_1$ to
\begin{eqnarray*}
    \left\{\begin{array}{rcl}
(\gamma u'_1)' &=&0\text{ in }(-1,1),  \\
 (\gamma u'_1) (-1)&=&(\gamma u'_1)(1)=K\in \mathbb C,\quad u_1(-1)+u_1(1)=0,
 \end{array}\right.
\end{eqnarray*}
one gets
\[
    u_1(x)=\left\{\begin{array}{ccl}

    F_0(x)+M&\mbox{if}& x\in (-1,a),\\
   F_1(x)+M+\frac{F_0(a)+F_0(b)}{2}-\frac{F_1(a)+F_1(b)}{2}&\mbox{if}&x\in (a,b),\\
  F_0(x)+M+\frac{F_1(b)-F_1(a)}{2}-\frac{F_0(b)-F_0(a)}{2}&\mbox{if}&x\in (b,1),
    \end{array}\right.\]
where
\[
F_1(x)=\int \frac{K}{\gamma_1(x)} dx.
\]
Hence
\begin{eqnarray*}
\dw&=&K\overline{(u_1(1)-u_0(1))-(u_1(-1)-u_0(-1))}=K\overline{(u_1(1)-u_0(1))}\\
&=&\frac{|K|^2}{2}\int_a^b\overline{\left(\frac{1}{\gamma_1}-\frac{1}{\gamma_0}\right)}dx,
\end{eqnarray*}
\[
\Re (\delta W)= \frac{|K|^2}{2}\int_a^b \left(\frac{\sigma_0}{\sigma_0^2+\epsilon_0^2}- \frac{\sigma_1}{\sigma_1^2+\epsilon_1^2}\right) dx,
\]
and
\[
\Im  (\delta W)=\frac{|K|^2}{2}\int_a^b \left(-\frac{\epsilon_0}{\sigma_0^2+\epsilon_0^2}+\frac{\epsilon_1}{\sigma_1^2+\epsilon_1^2}\right) dx.
\]
So, if one of the following monotonicity conditions hold
\[
 \frac{\sigma_1}{\sigma_1^2+\epsilon_1^2}>(<)\frac{\sigma_0}{\sigma_0^2+\epsilon_0^2} \text{ in } (-1,1)
\]
or
\[
 \frac{\epsilon_1}{\sigma_1^2+\epsilon_1^2}>(<)\frac{\epsilon_0}{\sigma_0^2+\epsilon_0^2}\text{ in} (-1,1)
\]
then either $\Re (\delta W) \neq 0 $  or   $\Im (\delta W)\neq 0$ and $\delta W$ recovers uniquely $(a,b)$.\\
In particular observe that if $\epsilon_0=0$ we find that $\Im (\delta W)\neq 0$ if $\epsilon_1$ has constant sign in$(-1,1)$ and $\Re (\delta W) \neq 0 $ if $\sigma_1-\sigma_0>0$ in $(-1,1)$ which are exactly conditions ${\bf (H3ii)}$.
Let us see that if the above conditions fail uniqueness does not hold. For, let us consider for example $\gamma_0=(2+ix)^2$ for $x\in (-1,1)$ and $\gamma_1=\frac{17}{4}$, then one easily sees that
\[
\Re (\delta W) =|K|^2(b-a)\left (\frac{4}{17}-\frac{4-ab}{(4+b^2)(4+a^2)}\right),
\]
\[
\Im (\delta W) =|K|^2(b-a)\left(-\frac{2(a+b)}{(4+b^2)(4+a^2)}\right)
\]
and clearly $\Re (\delta W)=\Im (\delta W)=0$ for $a=1/2$ and $b=-1/2$.

\section{Main tools: quantitative estimates of unique continuation}\label{sec4}
Let us list now various forms of the quantitative estimates of
unique continuation that we will need in the sequel.
Throughout this section we  will assume that $\om\subset \RR^n$ is a bounded domain of class $C^{0,1}$ with constants $r_0,M_0$.
and  $A$ is a symmetric $n\times n$
matrix in $\RR^n$ with real entries satisfying the following assumption

\vskip 1truemm

({\it Uniform ellipticity}) For a given $\lambda_0$, $0 <
\lambda_0 \leq 1$,
\begin{equation}
  \label{eq:2.ell}
  \lambda_0 | \xi | ^{2} \leq A(x) \xi \cdot \xi \leq {\lambda_0}^{-1} | \xi |
  ^{2}, \quad \quad \textrm{for every} \quad \xi \in\RR^{n}, x \in
  \RR^n.
\end{equation}

\vskip 1truemm

({\it Lipschitz regularity}) For a given $\lambda_1 > 0$
\begin{equation}
  \label{eq:2.Lip}
  | A(x) - A(y) | \leq \frac {\lambda_1} {r_0} | x - y |,
  \quad \textrm{for every} \quad x, y \in \RR^n.
\end{equation}
%

\begin{theo} [Three Spheres Inequality, \cite{AMR}]\label{TSI}
 Let $u\in H^1(\Omega)$
be a solution to the equation
\[
 \mathrm{div}(A(x)\nabla u(x))=0\mbox{ in }\om.
\]
For every
$r_1, r_2, r_3, \bar r$, $0<r_1<r_2<r_3\leq \bar r$, and for every
$x_0\in \Omega_{\bar r}$
\[
   \int_{B_{r_{2}}(x_0)}|{\nabla} u_0|^{2}
   \leq C
   \left(  \int_{B_{r_{1}}(x_0)}|{\nabla} u|^{2}
   \right)^{\theta}\left(  \int_{B_{r_{3}}(x_0)}|{\nabla} u|^{2}
   \right) ^{1-\theta},
\]
where $C>0$ and $\theta$, $0<\theta<1$, only depend on $\lambda_0$,
$\lambda_1$, $\frac{r_{1}}{r_{3}}$ and $\frac{r_{2}}{r_{3}}$.
\end{theo}

\begin{theo} [Lipschitz Propagation of Smallness, \cite{AMR}] \label{LPS}
\label{theo 4.2}  Let $h$ satisfy assumption {\bf (H4)} and let $u\in H^1(\Omega)$
be the solution of the Neumann problem
\begin{equation}
  \label{lipeq}
  \left\{ \begin{array}{ll}
   \mathrm{div}\left(A(x)\nabla u(x)\right)=0 &
  \mbox{in}\ \Omega ,\\
  &  \\
  A \nabla u \cdot \nu=h& \mbox{on}\ \partial
  \Omega.\\
  \end{array}\right.
\end{equation}
 For every
$\rho>0$ and for every $x\in\Omega_{2\rho}$, we have
\[
\int_{B_\rho(x)}|{\nabla} u|^2\geq
C^{-1}\int_\Omega |{\nabla} u|^2,
\]
where $C\geq 1$ only depends on $\lambda_0$, $\lambda_1$, $M_0$, $M_1$, $F(h)$ and $\frac{\rho}{r_0}$.
\end{theo}

The Three Spheres Inequality and  of the Lipschitz Propagation of Smallness in  \cite{AMR} are obtained for real valued functions $u$ and $h$ but with straightforward modifications they apply to complex valued functions.

\begin{theo} [Doubling Inequality]
\label{theo:4.dc} Let  $u\in H^1(B_{r_0}(x_0))$ be solution of
\begin{equation}
\label{doubeqeq}
\mathrm{div}(A(x)\nabla u(x))=0
 \quad \mathrm{ in }\ B_{r_0}(x_0).
\end{equation}
Then, there exist positive constants $\alpha, C$, depending only on $\lambda_0$ and on $\lambda_1$, such that
\begin{equation}
  \label{doubeq}
\int_{B_{2r}(x_0)}|{\nabla}u|^2 \leq C\left(\frac{\int_{B_{r_0}(x_0)}|{\nabla}u|^2}{\int_{B_{r_0/2}(x_0)}|{\nabla}u|^2}\right)^{\alpha}
\int_{B_r(x_0)}|{\nabla}u|^2,\end{equation}
for every $r$ such that $0<r\leq
\frac{ r_0} {2}$.

\end{theo}

The Doubling Inequality has been derived for the first time by Garofalo and Lin in \cite{GL}. Afterwards was also derived by Kukavica in \cite{Ku} using Rellich's identity.  In the Appendix, for convenience of the reader, we will give the proof of the Doubling Inequality following the proof on \cite{Ku}, showing the modifications one has to do in the case of complex valued functions and estimating more carefully the constant involved in the inequality.

\begin{lem}
\label{lemma} Let $\om$ satisfy assumption {\bf (H1)}, let $\overline{r}$ and $R$ be positive numbers such that $3\sqrt n R< \overline{r}$ and let $u\in H^1(\om)$ be a non trivial  solution of
\[
\mathrm{div}(A(x)\nabla u(x))=0
  \text{ in } \om.
\]
Assume that $\om_{\overline{r}}\neq \emptyset $. Then, for every $x_0\in\om_{\overline{r}}$ and for every measurable set $E$, $E\subset Q_R(x_0)$, we have

\begin{equation}
  \label{lemmaeq}
\frac{|E|}{|Q_R(x_0)|}\leq \left(\frac{H\int_E|\nabla u|^2}{\int_{Q_R(x_0)}|{\nabla}u|^2} \right)^{1/p},
\end{equation}
 where $H$ and $p>1$ are given by
 \[
 p=1+\frac{\log 4F_{\overline{r}}(u)}{\log (17/16)},
 \]
  \[
 H=\left(27 F_{\overline{r}}(u)\right)^{p(p-1)},
 \]
 where
\begin{equation}
\label{osc}
F_{\overline{r}}(u)=C\left(\frac{\int_{\om}|\nabla u|^2}{\int_{\om_{\frac{\overline{r}}{2}}}|\nabla u|^2}\right)^C
\end{equation}
and $C$ depends on $\lambda_0$, $\lambda_1$, $M_0$, $M_1$ and $\frac{\overline{r}}{r_0}$.
\end{lem}

\textit{A sketch of the proof.}
The proof of the above lemma can be derived with some changes in the proof of  \cite[Lemma 2.4]{Ve}. Let us illustrate the changes in detail. The most important difference between Lemma 2.4 in \cite{Ve} and our lemma is that in the bound (\ref{lemmaeq}) $|\nabla u|^2$ appears while in \cite{Ve} $|u|^2$ is involved.

Observe that  $|\nabla u|^2$   satisfies  the following reverse H\"{o}lder inequality (RHI)
\begin{equation}
\label{rhi}
\left(\frac{1}{|Q_R(x_0)|}\int_{Q_R(x_0)}\!\!\!\!\left(|\nabla u|^2\right)^{1+\delta}\right)^{\frac{1}{1+\delta}}\leq \frac{C}{|Q_R(x_0)|}\left(\frac{\int_{B_{\overline{r}}(x_0)} |\nabla u|^2}{\int_{B_{\overline{r}/2}(x_0)} |\nabla u|^2}\right)^{\alpha}\int_{Q_R(x_0)}|\nabla u|^2,
\end{equation}
for any $x_0\in \om_{\overline{r}}$ and $R$ such that $0<2\sqrt{n}R\leq \overline{r}$, where $C$ and $\alpha$ depend only on $\lambda_0$ and $\lambda_1$  and $\delta>0$ is arbitrary.

In fact, if we set $\tau= \frac{1}{|Q_R(x_0)|}\int_{Q_R(x_0)}u(x) dx$,  from \cite{GT}, Poincar\'e inequality and  (\ref{doubeq}) we get,
\begin{eqnarray*}
\sup_{Q_R(x_0)}|\nabla u|^2&\leq& \frac{C}{R^2}\sup_{Q_{3/2R}(x_0)}|u-\tau|^2\\
&\leq &\frac{C'}{R^{n+2}}\int_{Q_{2R}(x_0)}|u-\tau|^2  \\
&\leq& \frac{C''}{R^{n}}\int_{Q_{2R}(x_0)}|\nabla u|^2 \leq \frac{C''}{R^{n}}\int_{B_{2\sqrt{n}R}(x_0)}|\nabla u|^2 \\
&\leq& \frac{C'''}{R^{n}}\left(\frac{\int_{B_{\overline{r}}(x_0)}| \nabla u|^2}{\int_{B_{\overline{r}/2}(x_0)}|\nabla u|^2}\right)^{\alpha}\int_{Q_{R}(x_0)}|\nabla u|^2,
\end{eqnarray*}
where $C',C'',C'''$ and $\alpha$ depend on $\lambda_0,\lambda_1$ only. From last inequality we trivially derive (\ref{rhi}).\\
Using iteratively the Three Spheres Inequality we get the following estimate  (see \cite{ARRoV})
\begin{equation}
\label{properr}
\int_{\om_{\overline{r}/2}}|\nabla u|^2 \leq C \left(\int_{B_{\overline{r}/2}(x_0)}|\nabla u|^2 \right)^{\theta} \left(\int_{\om}|\nabla u|^2 \right)^{1-\theta},
\end{equation}
where $0<\theta<1$ and $\theta$ and $C$ depend on $\lambda_0,\lambda_1,M_0,M_1$ and $\frac{\overline{r}}{r_0}$.
From (\ref{properr}) we trivially have
\[
\frac{\int_{B_{\overline{r}}(x_0)}|\nabla u|^2}{\int_{B_{\overline{r}/2}(x_0)}|\nabla u|^2}\leq \frac{\int_{\om} |\nabla u|^2}{\int_{B_{\overline{r}/2}(x_0)} |\nabla u|^2}
\leq  \left(\frac{C\int_{\om} |\nabla u|^2}{\int_{\om_{\overline{r}/2}}|\nabla u|^2}\right)^{\frac{1}{\theta}}.
\]
From the above inequality and from (\ref{rhi}) we get the following version of RHI
\begin{equation}
\label{rhi'}
\left(\frac{1}{|Q_R(x_0)|}\int_{|Q_R(x_0)}\left(|\nabla u|^2\right)^{1+\delta}\right)^{\frac{1}{1+\delta}}\leq \frac{F}{|Q_R(x_0)|}\int_{Q_R(x_0)}|\nabla u|^2,
\end{equation}
for any $x_0\in \om_{\overline{r}}$ and for any $R$ such that $R\in (0,\frac{\overline{r}}{2\sqrt{n}}]$ and $\delta>0$, where
\begin{equation}
\label{defF}
F= \left(\frac{C\int_{\om} |\nabla u|^2 }{\int_{\om_{\frac{\overline{r}}{2}}}|\nabla u|^2}\right)^{\frac{\alpha}{\theta}},
\end{equation}
with $C,\alpha,\theta$ depending only  on $\lambda_0,\lambda_1,M_0,M_1$ and $\frac{\overline{r}}{r_0}$. In order to prove (\ref{properr}) we used the Lipschitz regularity of $\der\om$ in order to guarantee that $ \om_{\rho}$ is a connected set for $\rho$ sufficiently small. If (\ref{rhi'}) holds for $\overline{r}$ small then it clearly holds also for large $\overline{r}$ .
The most difficult part of the proof is to show that from (\ref{rhi'}) follows the thesis of the lemma but this part can be found in \cite{GC-RDeF} (Theorem 2.11) while an explicit evaluation of the constants can be found in \cite{Ve}.

\section{Proof of the main results}\label{sec5}

In this section we will use the quantitative unique continuation estimates stated in the previous section and regularity results for solutions of elliptic equations to get upper and lower bounds of the measure of the inclusion $D$, $|D|$, in terms of the energy related to  the background potential $u_0$. \\
Throughout this section we will assume that $A$ is a symmetric $n\times n$
matrix in $\RR^n$ satisfying (\ref{eq:2.ell}) and (\ref{eq:2.Lip})

We have
\begin{prop}\label{strictincl}
 Let $\om\subset \RR^n$ satisfy {\bf (H1)} with $\partial\om$ of class $C^{0,1}$, $D\subset\om$ satisfy {\bf (H2a)}.
Let $h$ satisfy {\bf (H4a)} and $u\in H^1(\om) $ solution to the Neumann problem (\ref{lipeq})
such that
\begin{equation}
\label{normeq}
\int_{\partial\om}u=0.
\end{equation}
Then
\begin{equation}
\label{uppbound}
\left(\frac{|D|}{|\om|}\right)^p\leq C\left(\frac{\int_D|\nabla u|^2}{\int_{\om}|\nabla u|^2}\right),
\end{equation}
where $p$, $p>1$, $C$ depend only on $\frac{d_0}{r_0}, M_0, M_1,\lambda_0,\lambda_1$ and on $F(h)$ where $F(h)$ is given by  (\ref{freq}).
 \end{prop}
\textit{Proof.}
Let $\delta=\frac{d_0}{4\sqrt n}$ and let us cover $D$ with closed cubes two by two internally disjoint, $Q_j$, $j=1,\cdots, N$ of size $\delta$. Assume that $Q_j\cap D\neq \emptyset$ for  $j=1,\ldots, N$. We have
\begin{equation}
\label{subset}
D\subset\cup_{j=1}^{N}Q_j\subset \om_{\frac{3}{4}d_0}.
\end{equation}

Let $p$, $p>1$, to be chosen later. From (\ref{subset}) and H\"{o}lder inequality (in what follows $p'=\frac{p-1}{p}$) we get
\begin{eqnarray*}
 |D|&=&\sum_{j=1}^N|D\cap Q_j|
=\sum_{j=1}^N\frac{|D\cap Q_j| }{|Q_j|}|Q_j|
\\
&\leq&\left(\sum_{j=1}^N\left(\frac{|D\cap Q_j| }{|Q_j|}\right)^p\right)^{1/p}\left(\sum_{j=1}^N|Q_j|^{p'}\right)^{1/p'}
\\
&\leq&|\om_{\frac{3}{4}d_0}|\left(\sum_{j=1}^N\left(\frac{|D\cap Q_j| }{|Q_j|}\right)^p\right)^{1/p}.
 \end{eqnarray*}
 Hence, for any $p>1$ we have
\begin{equation}
\label{uppeq1}
\left(\frac{|D|}{|\om|}\right)^p\leq |\sum_{j=1}^N\left(\frac{|D\cap Q_j| }{|Q_j|}\right)^p.
 \end{equation}
 Now, in order to choose $p$ and to bound the right-hand side of (\ref{uppeq1}) we apply Lemma 4.4 with $\overline{r}=\frac{3}{4}d_0$ and we bound from above $F_{\overline{r}}(h)$ defined in (\ref{osc}). For,  we bound from below $\int_{\om_{\frac{\overline{r}}{2}}}|\nabla u|^2$ observing that if $\bar x\in\om_{\overline{r}}$ then applying the Lipschitz Propagation of Smallness (LPS) with $\rho=\frac{\overline{r}}{2}$ we get
 \begin{equation}
 \label{lb1}
\int_{\om_{\frac{\overline{r}}{2}}}|\nabla u|^2\geq \int_{B_{\frac{\overline{r}}{2}}(\bar x)}|\nabla u|^2\geq C_1^{-1}\int_{\om}|\nabla u|^2,
\end{equation}
 where $C_1\geq 1$ depends on $\frac{d_0}{r_0}, M_0, M_1,\lambda_0,\lambda_1$ and on $F(h)$  where $F(h)$ is given by  (\ref{freq}). Hence, by (\ref{lb1}), we obtain that
 \[
 F_{\overline{r}}(h)\leq C_1.
 \]
 Let now
 \begin{equation}
 \label{eqp}
 p=1+\frac{\log 4C_1^2}{\log (17/16)}
 \end{equation}
 and by Lemma 4.4 we derive
 \begin{equation}
\label{uppbound1}
\left(\frac{|D\cap Q_j|}{|Q_j|}\right)^p\leq (27C_1^2)^{p(p-1)}\frac{\int_{D\cap Q_j}|\nabla u|^2}{\int_{Q_j}|\nabla u|^2},\quad j=1,\ldots, N.
\end{equation}
 Let us use  again the LPS property for estimating from above the right-hand side of (\ref{uppbound1}). Denoting by $x_j$ the center of the cube $Q_j$ we have
  \begin{equation}
 \label{lb2}
\int_{Q_j}|\nabla u|^2\geq \int_{B_{\frac{\delta}{2}}(x_j)}|\nabla u|^2\geq C_2^{-1}\int_{\om}|\nabla u|^2,
\end{equation}
 where $C_2\geq 1$ depends on $\frac{d_0}{r_0}, M_0, M_1,\lambda_0,\lambda_1$ and on $F(h)$. By (\ref{lb2}), (\ref{uppbound1}), (\ref{eqp}) and (\ref{uppeq1}) we get the thesis.
\qed
\begin{prop}\label{incl}
 Let $\om\subset \RR^n$ satisfy {\bf (H1)} with $\partial\om$ of class $C^{1,1}$, $D\subset\om$ satisfy {\bf (H2b)}, the function $\sigma_0$ as in {\bf (H3)} and  $h$ satisfying {\bf (H4b)}. Let  $u\in H^1(\om) $ solution to the Neumann problem
\[
  \left\{ \begin{array}{ll}
   \mathrm{div}((\sigma_0(x)\nabla u(x))=0 &
  \mathrm{in}\ \Omega ,\\
  &  \\
  \sigma_0 \nabla u \cdot \nu=h& \mathrm{on}\ \partial
  \Omega,\\
  \end{array}\right.
\]
satisfying the normalization condition (\ref{normeq}).
Then
\[
\left(\frac{|D|}{|\om|}\right)^p\leq C\left(\frac{\int_D|\nabla u|^2}{\int_{\om}|\nabla u|^2}\right),
\]
where $p$, $p>1$, $C$ depend only on $\frac{r_1}{r_0}, M_0, M_1,c_0, L$ and on $F(h)$ where $F(h)$ is given by  (\ref{freq}).
 \end{prop}

\textit{Proof.}
Denote by  $\Gamma:=\partial\Omega\cap\overline{B_{r_1}(P)}$.
First we construct a suitable family of cylinders covering $\partial\om\backslash\Gamma$.\\
Let
\[
r_2=\min\left\{\frac{r_1}{4\sqrt n},\frac{r_0}{2\sqrt n M_0}\right\}
\]
and let $r\in (0,r_2]$ to be chosen later. Let $\{Q_j\}_{j=1}^J$ a family of closed mutually internally disjoint cubes of size $2r$ such that
\[
 (\partial\om\backslash\Gamma)\cap Q_j\neq\emptyset,\quad j=1,\ldots, J,
\]
\[
\partial\om\backslash\Gamma\subset \cup_{j=1}^J Q_j.
\]
Let us fix $j\in \{1,\ldots, J\}$ and let $x_j\in (\partial\om\backslash\Gamma)\cap Q_j$. Let $\nu_j$ the exterior normal unit vector in $x_j$ on $\partial\om$. Let $R_j$ the cylinder centered at $x_j$ with axis parallel to $\nu_j$ and with basis a ball of radius equal to $2\sqrt n r$ and with height equal to $2\sqrt n M_0 r$. Setting $\tilde{R_j}=2(R_j-x_j)+x_j$ one easily see that
\[
\cup_{j=1}^J \tilde{R_j}\supset \om\backslash\om_{2\sqrt nr}
\]
and hence
\begin{equation}
\label{cyl2}
\textrm{dist}(\om\backslash\cup_{j=1}^J \tilde{R_j},\partial\om) \geq 2\sqrt n r.
\end{equation}
Furthermore, since the interior of the cubes $Q_j$,  $ j=1,\ldots, J$, are two by two disjoint and since, obviously
\[
\cup_{j=1}^J Q_j\subset \om\backslash\om_{2\sqrt nr},
\]
we derive the following estimate for $J$
\begin{equation}
\label{cyl4}
J \leq (2r)^{-n}|\cup_{j=1}^J Q_j|\leq C\left(\frac{r_0}{r}\right)^{n-1},
\end{equation}
where $C$ depends only on $M_0,M_1$. Let

\[
D'=D\cap (U_{j=1}^J\tilde{R_j}),
\]
\[
D''=D\backslash D'.
\]
From (\ref{cyl2}) we have
\[
\textrm{dist}(D'',\der\om)\geq 2\sqrt nr.
\]
From last inequality and Proposition 5.1 we get
\begin{equation}
\label{uppbound2}
\left(\frac{|D''|}{|\om|}\right)^p\leq C_r\left(\frac{\int_{D''}|\nabla u|^2}{\int_{\om}|\nabla u|^2}\right),
\end{equation}
where $C_r$ depends only on $\frac{r}{r_0}, M_0, M_1,c_0, L$ and on $F(h)$.\\
Let, for a fixed index $j\in \{1,\ldots, J\}$ ,
\[
D_j:=\tilde{R_j}\cap D'
\]
and
\[
\hat{R}_j:=2(\tilde{R_j}-x_j)+x_j.
\]
It is easy to see that if
\[
r\leq \min\left\{\frac{r_1}{16\sqrt n\sqrt{1+M_0^2}}\,,\frac{r_2}{2}\right\}
\]
then
\begin{equation}
\label{cyl5}
\textrm{dist}(\hat{R}_j,\Gamma_0)\geq\frac{r_1}{4},
\end{equation}
where, we recall that, $\Gamma_0=\Gamma\cap \overline{B_{r_1/2}(P)}$. Furthermore, up to a rigid transformation such that with $x_j=0$, we have
\[
\hat{R}_j\cap\om=\{(x',x_n)\in\RR^n: x_n>\psi (x'),|x'|\leq 8\sqrt n r, |x_n|\leq 8\sqrt n M_0 r\},
\]
where
\[
\psi(0)=|\nabla\psi(0)|=0
\]
and
\[
\|\psi\|_{L^{\infty}}+r_0\|\nabla\psi\|_{L^{\infty}}+r_0^2\|D^2\psi\|_{L^{\infty}}\leq M_0r_0.
\]
Without loss of generality we may assume that $\sigma_0(0)=1$. Following the arguments of \cite{AE} or  \cite{ABRV} we can construct a function $\Psi\in C^{1,1}(\overline{B_{\rho_0}(0)},\RR^n)$, where $\rho_0=16\sqrt n\sqrt{1+M_0^2}r$ such that
\begin{equation}
\label{Psi1}
\Psi(x',\psi(x'))=(x',0),\quad \forall x'\in B'_{\rho_0}(0),
\end{equation}
\begin{equation}
\label{Psi2}
\Psi(\hat{R}_j\cap\om)\subset \{(x',x_n)|x_n>0\}.
\end{equation}
Moreover, there exist $C_1,C_2\geq 1$ depending only on $M_0$ such that
\begin{equation}
\label{Psi3}
C_1^{-1} |x-z|\leq |\Psi(x)-\Psi(z)|\leq C_1|x-z|,\quad \forall x,z\in B_{\rho_0}(0),
\end{equation}
\begin{equation}
\label{Psi4}
C_2^{-1}\leq |\textrm{det}D\Psi(x)|\leq C_2,\quad \forall x\in B_{\rho_0}(0),
\end{equation}
and, setting $A(y)=\{a_{ij}(y)\}_{ij=1}^n$, where
\begin{equation}
\label{Psi5}
A(y)=|\textrm{det}D\Psi^{-1}(x)|(D\Psi)(\Psi^{-1}(y))\sigma_0(\Psi^{-1}(y))(D\Psi)^T(\Psi^{-1}(y)),
\end{equation}
\begin{equation}
\label{Psi6}
v(y)=u(\Psi^{-1}(y)),
\end{equation}
we have
\begin{equation}
\label{A1}
A(y)=Id,
\end{equation}
\begin{equation}
\label{A2}
a_{nk}(y',0)=a_{kn}(y',0)=0, k=1,\ldots,n,
\end{equation}
\begin{equation}
\label{A3}
C_3^{-1}|\xi|^2\leq A(y)\xi\cdot \xi\leq C_3|\xi|^2,\quad\forall\xi\in\RR^n,\forall y\in\Psi(\om\cap\hat{R}_j),
\end{equation}
\begin{equation}
\label{A4}
|A(y)-A(z)|\leq \frac{C_4}{r}|y-z|,\quad\forall y,z\in\Psi(\om\cap\hat{R}_j),
\end{equation}
where in (\ref{A1}) $Id$ denotes the identity matrix and $C_3,C_4\geq 1$ depend only on $M_0$. Furthermore, recalling that $\hat{R}_j\cap\Gamma_0\neq\emptyset$ and (\ref{cyl5}) we have
\begin{equation}
  \label{v1}
  \left\{ \begin{array}{ll}
   \mathrm{div}((A(y)\nabla_y v(y))=0 &
  \mathrm{in }\Psi(\om\cap\hat{R}_j) ,\\
  &  \\
\frac{\partial v}{\partial y_n}(y',0) =0& \mathrm{on}\ \ \Psi(\der\om\cap\hat{R}_j).\\
  \end{array}\right.
\end{equation}
From the properties of the matrix $A$, in particular from (\ref{A2}), we have that the function $\tilde v$ defined by
\begin{equation}
\label{v2}
\tilde v(y',y_n):=v(y',|y_n|)
\end{equation}
is solution of an elliptic equation with Lipschitz coefficients in the principal part. More precisely, let $\tilde A(y)=\{\tilde a_{ij}(y)\}_{i,j=1}^n$ the matrix with entries given by
\[
\tilde a_{ij}(y',|y_n|)=a_{ij}(y',|y_n|), \textrm{if }  i,j\in\{1,\ldots,n-1\}\textrm{ or }i=j=n,
\]
\[
\tilde a_{ij}(y',y_n)=\tilde a_{ij}(y',y_n)=\textrm{sgn}(y_n)a_{nj}(y',|y_n|),\textrm{if }  i,j\in\{1,\ldots,n-1\}\textrm{ or }i=j
\]
then we have
\[
  \mathrm{div}((A(y)\nabla_y \tilde v(y))=0  \textrm{ in }\hat{\Lambda}_j,
\]
where
\[
\hat{\Lambda}_j=\{(y',y_n)\in\RR^n|\,(y',|y_n|)\in \hat{\Lambda}_j^+\},
\]
with
\[
\hat{\Lambda}_j^+=\Psi(\om\cap\hat{R}_j).
\]
It is easy to see that the matrix $\tilde A$ satisfies uniform ellipticity and Lipschitz continuity with the same constants as in (\ref{A3}) and (\ref{A4}). \\
In the sequel we will use the following notations
\[
\tilde{\Lambda}_j^+:=\Psi(\om\cap\tilde{R}_j),
\]
\[
\tilde{\Lambda}_j:=\{(y',y_n)\in\RR^n|\,(y',|y_n|)\in \tilde{\Lambda}_j^+\},
\]
\[
\tilde{D}_j:=\Psi(D_j).
\]
Since our aim is to bound $|D_j|$, let us proceed initially as in the proof of Proposition 5.1.\\
Let us note first that from (\ref{Psi3}) we get
\[
\textrm{dist}(\tilde{D}_j,\der\hat{\Lambda}_j)\geq \delta_0:=\frac{2\sqrt n r}{C_1},
\]
where $C_1$ is the constant appearing in (\ref{Psi3}). Let us cover $\tilde{D}_j$ with closed cubes, two by two internally disjoint, $Q_{j,k}, k=1,\ldots, N_j$ of size $\delta_1:=\frac{\delta_0}{4\sqrt n}$. We have
\[
\tilde{D}_j\subset \cup_{k=1}^{N_j}Q_{j,k}\subset \hat{\Lambda}_{\frac{3}{4}\delta_0}.
\]
Since we are interested in applying Lemma 4.4 with $\om=\hat{\Lambda}_j$ and $\bar r=\frac{\delta_0}{4}$ we need to prove first the following
\begin{claim}
There exists a constant $C$ depending only on $c_0, L,M_0,M_1,\frac{r}{r_0}$ and $F(h)$ such that
\begin{equation}
\label{tildeF1}
\tilde{F}_{j,\bar r}(\tilde v):=\frac{\int_{\hat{\Lambda}_j}|\nabla\tilde v|^2}{\int_{\hat{\Lambda}_{j,\bar r/2}}|\nabla\tilde v|^2}\leq C, \quad j=1,\ldots,J,
\end{equation}
(with $C$ independent on $j$).
\end{claim}
\textit{Proof of the claim.}

Since, for $\bar r=\delta_0/4$ we have that $\hat{\Lambda}_{j, \frac{\bar r}{2}}\supset \tilde{\Lambda}$, recalling that $\tilde v$ is the even reflection of $v=u\circ\Psi^{-1}$, by a change of variables we derive
\begin{equation}
\label{tildeF2}
\tilde{F}_{j,\bar r}(\tilde v)\leq C\frac{\int_{\hat{R}_j\cap\om}|\nabla  u|^2}{\int_{\tilde{R}_{j,\bar r/2}\cap\om}|\nabla u|^2},\quad j=1,\ldots,J,
\end{equation}
where $C$ depends only on $c_0, L,M_0,M_1,\frac{r}{r_0}$. Now, since
\[
\tilde{R}_j\cap\om\supset B_{\sqrt n r}(x_j-2\sqrt n r\nu):=B^{(j)},
\]
\[
\textrm{dist} (B^{(j)}, \der(\tilde{R}_j\cap\om))\geq \sqrt n r,
\]
estimating the right-hand side of (\ref{tildeF2}) and applying the LPS property we get
\[
\tilde{F}_{j,\bar r}(\tilde v)\leq C\frac{\int_{\om}|\nabla u|^2}{\int_{B^{(j)}}|\nabla u|^2}\leq C',
\]
where $C'$ depends only on $c_0, L,M_0,M_1,\frac{r}{r_0}$ and $F(h)$.\\

Let us choose $r=r_2$. Proceeding  as in the proof of Proposition 5.1 and using (\ref{tildeF2}) we obtain
\begin{equation}
\label{tildeDj}
|\tilde{D}_j|\leq |\hat{\Lambda}_j|\left(\frac{\int_{\tilde{D}_j}|\nabla\tilde v|^2}{\int_{\hat{\Lambda}_{j}}|\nabla\tilde v|^2}\right)^{1/p}, \quad j=1,\ldots,J,
\end{equation}
where $C$ and $p\in (1,+\infty)$ depend on $c_0, L,M_0,M_1,\frac{r_1}{r_0}$ and $F(h)$.\\
From the definition of $\tilde v$ and of $\hat{\Lambda}_j$, with some simple change of variables and using again the LPS property, we derive from (\ref{tildeDj})
\begin{equation}
\label{Dj}
|D_j|\leq C|\om|\left(\frac{\int_{D_j}|\nabla u|^2}{\int_{\om}|\nabla u|^2}\right)^{1/p}, \quad j=1,\ldots,J,
\end{equation}
where $C$ and $p\in (1,+\infty)$ depend on $c_0, L,M_0,M_1,\frac{r_1}{r_0}$ and $F(h)$.

From (\ref{Dj}) and from (\ref{cyl4}) we have
\begin{equation}
\label{D'}
|D'|\leq \sum_{j=1}^J|D_j|\leq C|\om|\left(\frac{\int_{D'}|\nabla u|^2}{\int_{\om}|\nabla u|^2}\right)^{1/p},
\end{equation}
where $C$ and $p\in (1,+\infty)$ depend on $c_0, L,M_0,M_1,\frac{r_1}{r_0}$ and $F(h)$. From (\ref{D'}) and (\ref{uppbound2}) the thesis follows.
\qed
\begin{prop}\label{size}
Under the same hypotheses of Proposition 5.2 we have
\[
\frac{\int_D|\nabla u|^2}{\int_{\om}|\nabla u|^2}\leq Cr_1^{-n}|D|
\]
where $C$ depends on $\frac{r_1}{r_0}, M_0, M_1,c_0, L$ only.
 \end{prop}
\textit{Proof.}
We trivially have
\begin{equation}
\label{energy1}
\int_D|\nabla u|^2\leq |D|\,\|\nabla u\|^2_{L^{\infty}(D)}\leq |D|\,\|\nabla u\|^2_{L^{\infty}(\om\backslash B_{r_1}(P))}.
\end{equation}
Since $\sigma_0\nabla u\cdot \nu=0$ on $\der\om\backslash  B_{r_1/2}(P)$,  from standard estimates for elliptic equations, \cite{GT}, and from Poincar\'e inequality, we have that, letting $\tau=\frac{1}{|\om|}\int_{\om}u$, the following inequalities hold
\begin{eqnarray}
\label{ineq}
\|\nabla u\|^2_{L^{\infty}(\om\backslash B_{r_1}(P))}&\leq&\frac{C_1}{r_1^2}\|u-\tau\|^2_{L^{\infty}(\om\backslash B_{\frac{3}{4}r_1}(P))}\nonumber\\
&\leq&\frac{C_1C_2}{r_1^{n+2}}\|u-\tau\|^2_{L^2(\om\backslash B_{r_{1/2}}(P))}\\
&\leq&\frac{C_1C_2}{r_1^{n+2}}\|u-\tau\|^2_{L^2(\om)}\leq \frac{C_1C_2C_3 r_0^2}{r_1^{n+2}}\|\nabla u\|^2_{L^2(\om)},\nonumber
\end{eqnarray}
where $C_1$ depends only on $\frac{r_1}{r_0}, M_0, M_1,c_0, L$, $C_2$ depends on $\frac{r_1}{r_0}, M_0, M_1,c_0$ and $C_3$ depends on $ M_0, M_1$. From (\ref{energy1}) and (\ref{ineq}) we get
\[
\frac{\int_D|\nabla u|^2}{\int_{\om}|\nabla u|^2}\leq C_4 r_1^{-n}|D|,
\]
where $C_4$ depends only on $\frac{r_1}{r_0}, M_0, M_1,c_0, L$.\qed

We are now ready to prove our main results\\

\textit{Proof of Theorem 2.1.}
By standard elliptic estimates, we have
\[
\sup_D|\nabla u_0|\leq C\sup_{\Omega_{d_0 / 2}} |u_0|\leq C
\|u_0\|_{L^2(\Omega)},
\]
From the trivial estimate
\[
\int_{\Omega} |\nabla u_0|^2\leq c_0^{-1} \int_{\Omega} \nabla u_0\cdot \nabla u_0= c_0^{-1}W_0,
\]
and from Poincar\'e inequality, we have
\begin{equation}
  \label{triv2}
\sup_D|\nabla u_0|\leq CW_0^\frac{1}{2},
\end{equation}
where $C$ only depends only on $c_0$, $L$, $\frac{d_0}{r_0}$ and $M_0$. Hence from (\ref{triv2}) we get a lower bound for $|D|$
\begin{equation}
\label{energy1'}
\int_D|\nabla u|^2\leq |D|\,\|\nabla u\|^2_{L^{\infty}(D)}\leq C |D|W_0.
\end{equation}
By Proposition 5.1 and (\ref{energy1'}) we obtain
\[
\tilde{C}_1\frac{\int_D|\nabla u_0|^2}{W_0}\leq\frac{|D|}{|\om|}\leq \tilde{C}_2\left(\frac{\int_D|\nabla u_0|^2}{W_0}\right)^{1/p},
\]
where $\tilde{C}_1$ depends only on $\frac{d_0}{r_0}$, $M_0$, $M_1$, $c_0$, $L$ and $\tilde{C}_2$ depends only on $\frac{d_0}{r_0}$, $M_0$, $M_1$, $c_0$, $L$ and on $F(h)$. Finally, applying Proposition 3.1,  if $\gamma_0,\gamma_1$ are constant and satisfy  {(\bf H3i)} or Proposition 3.2 if $\gamma_0$ and $\gamma_1$ satisfy {(\bf H3ii)} we get
\[
C_1\left|\frac{\delta W}{W_0}\right|\leq \frac{|D|}{|\om|}\leq C_2\left|\frac{\delta W}{W_0}\right|^{1/p},
\]
where $C_1$ depends on the a priori constants $c_0$, $\mu_0$, $M_0$, $M_1$, $\frac{d_0}{r_0}$, $L$, and the numbers $p>1$ only and $C_2$ depends on the same a priori parameters and, on $F(h)$ only.
\qed

\textit{Proof of Theorem 2.2.}
By Proposition 5.2 and Proposition 5.3 applied to the background potential $u_0$ solution of (\ref{5}) (in the constant case up to a rescaling by a constant) we get
\[
C'_1\frac{\int_D|\nabla u_0|^2}{W_0}\leq\frac{|D|}{|\om|}\leq C'_2\left(\frac{\int_D|\nabla u_0|^2}{W_0}\right)^{1/p},
\]
where $C'_1$ depends only on $\frac{r_1}{r_0}$, $M_0$, $M_1$, $c_0$, $L$ and $C'_2$ depends on $\frac{r_1}{r_0}$, $M_0$, $M_1$, $c_0$, $L$ and on $F(h)$ only. Finally applying Proposition 3.1 if $\gamma_0,\gamma_1$ are constant and satisfy  {(\bf H3i)} or Proposition 3.2 if $\gamma_0$ and $\gamma_1$ satisfy {(\bf H3ii)} we get
\[
C_1\left|\frac{\delta W}{W_0}\right|\leq \frac{|D|}{|\om|}\leq C_2\left|\frac{\delta W}{W_0}\right|^{1/p},
\]
where $C_1$ depends only on the a priori constants $c_0, \mu_0, M_0, M_1, \frac{r_1}{r_0}, L$, and the numbers $p>1$ only and $C_2$ depends only on the same a priori parameters and, on $F(h)$.
\qed

\section{Appendix}\label{sec7}
\textit{Proof of Proposition 4.1.}
The Doubling Inequality proved in \cite{Ku} and  \cite{ARRoV} can be extended with straightforward arguments to the case of complex valued solutions of
\[
\mbox{div}(A(x)\nabla u(x))=0 \quad
  \mathrm{ in }\ B_{r_0}(x_0).
\]
let us give an idea of the modifications to be done in the proof.  Let us assume that
\begin{equation}
\label{norm1}
A(0)=Id
\end{equation}
and let, for $0<r<R_0$,
\begin{equation}
\label{def}
H(r)=\int_{\der B_r}\frac{A(x)x\cdot x}{|x|^2}|v(x)|^2,\quad
I(r)=\int_{B_r}A(x)\nabla v\cdot\overline{\nabla v},\quad
N(r)=\frac{rI(r)}{H(r)}.
\end{equation}
If, instead of Rellich's identity used in \cite{Ku} we use the following relation
\begin{eqnarray*}
2\Re\left[\vphantom{\int_L}(\beta\cdot\nabla\bar v)\textrm{div}(A\nabla v)\right]&=&
\textrm{div}\left[\vphantom{\int_L}2\Re (( \beta\cdot\nabla\bar v)A\nabla v)-\beta (A\nabla v\cdot\nabla\bar v)\right]\\
+(\textrm{div}\beta)A\nabla v\cdot\nabla\bar v&-&2\Re\left[\vphantom{\int_L}\der_l\beta_ja_{lk}\der_k v\der_j\bar v\right]
+\beta_j(\der_ja_{lk})\der_kv\der_l\bar v,
\end{eqnarray*}
with $\beta$ vector field in $\RR^n$ sufficiently smooth, we get that there exist constant $C_1>1,C_2$ and $c$, with $C_1>1,C_2$ depending only on $\lambda_0,\lambda_1$ and $c$ absolute constant, such that
\begin{equation}
\label{ineqA1}
\int_{B_r}|v|^2\leq \lambda_0^2 r\int_{\der B_r}|v|^2, \textrm{ for }r\leq \frac{R_0}{C_1},
\end{equation}
\begin{equation}
\label{ineqA2}
|H'(r)-\frac{n-1}{r}H(r)-2I(r)|\leq \frac{c\lambda_1}{R_0}H(r),
\end{equation}
\begin{equation}
\label{ineqA3}
N(r)e^{C_2\frac{r}{R_0}}\quad\textrm{increasing in }(0,R_0].
\end{equation}
From (\ref{ineqA2}) we have
\begin{equation}
\label{ineqA5}
\frac{d}{dr}\left(\log \frac{H(r)}{r^{n-1}}\right)\leq \frac{c\lambda_1}{R_0}+\frac{2N(r)}{r},
\end{equation}
\begin{equation}
\label{ineqA6}
\frac{2N(r)}{r}\leq \frac{d}{dr}\left(\log \frac{H(r)}{r^{n-1}}\right)+\frac{c\lambda_1}{R_0}.
\end{equation}
Let $R_1:=\frac{R_0}{C_1}$ and $\rho,R\in (0,R_1]$ such that $3\rho\leq R$. Integrating both sides of (\ref{ineqA5}) in the interval $[\rho,3\rho]$ we get, using (\ref{ineqA3}),
\begin{eqnarray*}
\log \frac{H(3\rho)}{3^{n-1}H(\rho)}&\leq& \frac{2c\lambda_1\rho}{R_0}+\int_{\rho}^{3\rho}\frac{2N(r)}{r}\\
&\leq& \frac{2c\lambda_1\rho}{R_0}+\int_{\rho}^{3\rho}\frac{2N(r)}{r} e^{C_2\frac{r}{R_0}}\\
&\leq& \frac{2c\lambda_1\rho}{R_0}+2N(3\rho)e^{3C_2\frac{\rho}{R_0}}\log 3\\
&\leq& \frac{2c\lambda_1R}{3R_0}+2N(R)e^{3C_2\frac{R}{R_0}}\log 3.
\end{eqnarray*}
Hence, for $\rho\in(0,R/3]$, $R\in (0,R_1]$, one has
\begin{equation}
\label{ineqA7}
\frac{1}{R}\log \frac{H(3\rho)}{3^{n-1}H(\rho)}\leq \frac{2c\lambda_1}{3R_0}+2e^{C_2\frac{R}{R_0}}\frac{N(R)}{R}\log 3.
\end{equation}
From (\ref{ineqA7}) and (\ref{ineqA6})  we get, for $\rho\in(0,R_1/3]$, $R\in (0,R_1]$,
\begin{equation}
\label{ineqA8}
\frac{1}{R}\log \frac{H(3\rho)}{3^{n-1}H(\rho)}\leq \frac{C_3}{R_0}+e^{C_2}(\log 3)\frac{d}{dR}\left(\log \frac{H(R)}{R^{n-1}}\right),
\end{equation}
where $C_3=e^{C_2\log 3}+\frac{2c\lambda_1}{3}$. Last inequality implies, in particular, that for any  $\rho\in(0,R_1/9]$, $R\in (R_1/2,3R_1/4]$ one has (integrating both sides of (\ref{ineqA8}) on $[R_1/3,R]$)
\begin{eqnarray*}
\log 6\log \frac{H(3\rho)}{3^{n-1}H(\rho)}&\leq& \log \frac{R}{R_1/3}\int_{R_1/3}^R\frac{1}{t}\log \frac{H(3\rho)}{3^{n-1}H(\rho)} dt\\
&\leq& C_3\lambda_1 \frac{R-R}{R_0}+e^{C_2}(\log 3)\log \frac{H(R)}{(\frac{3R}{R_1})^{n-1}H(R_1/3)}\\
&\leq&C_3\lambda_1+e^{C_2}(\log 3)\log \frac{H(R)}{(\frac{3}{2})^{n-1}H(R_1/3)}.
\end{eqnarray*}
Hence, for $\rho\in(0,R_1/9]$, $R\in (R_1/2,3R_1/4]$, by the elementary properties of the logarithmic function, we have,
\begin{equation}
\label{ineqA9}
H(3\rho)\leq C_4\left( \frac{H(R)}{H(R_1/3)}\right)^{C_5}H(\rho),
\end{equation}
where $C_4,C_5$ depend only on $\lambda_0,\lambda_1$. Integrating both sides of (\ref{ineqA9}), we derive, for every $\rho\in(0,R_1/9]$, $R\in (R_1/2,3R_1/4]$,
\[
\int_0^{\rho}H(3s)ds\leq C_4\left( \frac{H(R)}{H(R_1/3)}\right)^{C_5}\int_0^{\rho}H(s)ds.
\]
From (\ref{def}) we get
\[
\int_0^{\rho}H(s) ds\leq \lambda_0^{-1}\int_{B_{\rho}}|v|^2
\]
and
\[
\int_0^{\rho}H(3s) ds\geq \frac{\lambda_0}{3}\int_{B_{3\rho}}|v|^2.
\]
From the last two  inequalities and from (\ref{ineqA9}) one has
\begin{equation}
\label{ineqA11}
\int_{B_{3\rho}}|v|^2\leq 3\lambda_0^{-2}C_4\left( \frac{H(R)}{H(R_1/3)}\right)^{C_5}\int_{B_{\rho}}|v|^2,
\end{equation}
for $\rho\in(0,R_1/9]$, $R\in (R_1/2,3R_1/4]$.\\
Now, (\ref{ineqA11}) holds also if instead of $v$ we insert $v-\tau_{\rho}$ being $\tau_{\rho}=\frac{1}{|B_{\rho}|}\int_{B_{\rho}}v$. Denoting by $\tilde{H}(r)$ the function obtained substituting
in  (\ref{def}) $v-\tau_{\rho}$ instead of $v$ we have, recalling the theorem of local boundness of solutions of elliptic equations,  \cite{GT},
\begin{equation}
\label{ineqA12}
\tilde{H}(R)\leq\lambda_0^{-2}\int_{\der B_R}|v-\tau_{\rho}|^2\leq 4\lambda_0^{-2} R^{n-1}\|v\|_{L^{\infty}(B_R)}\leq
 \frac{C}{R_1}\int_{B_{R_1}}|v|^2,
\end{equation}
for every $R\in (R_1/2,3R_1/4]$ and $C$ depends only on $\lambda_0$. On the other hand, applying (\ref{ineqA1}), we derive
\begin{equation}
\label{ineqA13}
\tilde{H}(R_1/3)\geq\lambda_0\int_{\der B_{R_1/3 }}|v-\tau_{\rho}|^2
\geq \int_{B_{R_1/3}}|v-\tau_{\rho}|^2\geq
 \frac{R_1}{C}\int_{B_{R_1/6}}|\nabla v|^2,
\end{equation}
where $C\geq 1$ depends only on $\lambda_0$.
From what obtained in (\ref{ineqA12}),  (\ref{ineqA13}) and (\ref{ineqA11}) we get
\begin{equation}
\label{ineqA14}
\int_{B_{3\rho}}|v-\tau_{\rho}|^2
\leq C_6\left(\frac{\int_{B_{R_1}}|v|^2}{R_1^2\int_{B_{R_1/6}}|v|^2}\right)^{C_5}\int_{B_{\rho }}|v-\tau_{\rho}|^2,
\end{equation}
where $C_6\geq 1$ depends on $\lambda_0,\lambda_1$. Using Poincar\'e inequality and Caccioppoli inequality to bound from above the right hand side of (\ref{ineqA14}) and to bound from below the left hand side of (\ref{ineqA14}), we obtain for any $\rho\in(0,R_1/3]$,
\begin{equation}
\label{ineqA15}
\int_{B_{2\rho}}|\nabla v|^2
\leq C_7\left(\frac{\int_{B_{R_1}}|v|^2}{R_1^2\int_{B_{R_1/6}}|\nabla v|^2}\right)^{C_5}\int_{B_{\rho }}|\nabla v|^2,
\end{equation}
where $C_6\geq 1$ depends on $\lambda_0,\lambda_1$.\\
Iterating (\ref{ineqA15}) and by simple calculations we get
\begin{equation}
\label{ineqA16}
\int_{B_{\alpha\rho}}|\nabla v|^2
\leq C_8N_v'\alpha^{\frac{\log N_v}{\log 2}}\int_{B_{\rho }}|\nabla v|^2,
\end{equation}
for any $\alpha\geq 1$ and $\rho$ such that $3\alpha\rho\leq R_1$. Here we have set
\[
N'_v=\left(\frac{\int_{B_{R_1}}|v|^2}{R_1^2\int_{B_{\frac{R_1}{6}}}|\nabla v|^2}\right)^{C_5}
\]
and $C_8$ depends on $\lambda_0,\lambda_1$ only. Now let us remove condition (\ref{norm1}). Hence, let $A(x)$ be a symmetric matrix satisfying (\ref{eq:2.ell}) and (\ref{eq:2.Lip}) and let $v\in H^1(B_{R_0})$ a weak solution of equation (\ref{doubeqeq}). Let us introduce the change of variables $y=Jx$ where $J=\sqrt{A^{-1}(0)}$ and consider, for any $r>0$, the ellipsoids
\[
E_r:= \{x\in \RR^n | A^{-1}(0)x\cdot x<r^2\}=J^{-1}(B_r).
\]
Setting $w(y)=v(J^{-1}y)$ and $\tilde{A}(y)=JA(J^{-1}y)J$ one has
\begin{eqnarray*}
&\textrm{div}\left(\tilde{A}(y)\nabla_y w(y)\right)=0\textit{ in }B_{R_0\sqrt{\lambda_0}},&\\
&\lambda_0^2|\xi|^2\leq \tilde{A}(y)\xi\cdot\xi\leq\lambda_0^{-2}|\xi|^2, \forall y\in\RR^n,\forall\xi\in\RR^n,&\\
&|\tilde{A}(y_1)-\tilde{A}(y_2)|\leq \frac{\lambda_0^{-3/2}\lambda_1}{R_0}|y_1-y_2|,\forall y_1,y_2\in \RR^n \textrm{ and }\tilde{A}(0)=Id.&
\end{eqnarray*}
Furthermore, since
\[
B_{\sqrt{\lambda_0}r}\subset E_r\subset B_{\frac{r}{\sqrt{\lambda_0}}},\quad \forall r>0,
\]
by simple changes of variables we have
\begin{equation}
\label{ineqA17}
\lambda_0^{n/2+1}\int_{B_{\sqrt{\lambda_0}r}}|\nabla v|^2 dx\leq \int_{B_r}|\nabla w|^2 dy\leq \lambda_0^{-(n/2+1)}\int_{B_{\frac{r}{\sqrt{\lambda_0}}}}|\nabla v|^2 dx .
\end{equation}
for all $r>0$. We can apply (\ref{ineqA16}) to $w$ and we have
\begin{equation}
\label{ineqA18}
\int_{B_{\alpha\rho}}|\nabla w|^2 dy\leq C_8' N''_w\alpha^{\frac{\log N''_w}{\log 2}} \int_{B_{\rho}}|\nabla w|^2 dy
\end{equation}
for any $\alpha\geq 1$ and $\rho$ such that $3\alpha\rho\leq R_1\sqrt{\lambda_0}:=R_2$, where
\[
N''_w=\left(\frac{\int_{B_{R_2}}|w|^2dy}{R_2^2\int_{B_{R_2/6}}|\nabla w|^2}\right)^{C'_5}
\]
and $C'_5, C'_8$ depend on $\lambda_0,\lambda_1$ only. From (\ref{ineqA17}) and (\ref{ineqA18}) we derive easily
\begin{eqnarray}
\label{ineqA20}
\int_{B_{\alpha\rho}}|\nabla v|^2 dx&\leq& \lambda_0^{-(n/2+1)}\int_{B_{\frac{\alpha\rho}{\sqrt{\lambda_0}}}}|\nabla w|^2 dy\nonumber\\
&\leq& C'_8 N''_w(\lambda_0^{-1}\alpha)^{\frac{\log N''_w}{\log 2}}\lambda_0^{-(n/2+1)}\int_{B_{\rho\sqrt{\lambda_0}}}|\nabla w|^2 dy \nonumber\\
&\leq& C'_8\lambda_0^{-(n+2)}N''_w(\lambda_0^{-1}\alpha)^{\frac{\log N''_w}{\log 2}} \int_{B_{\rho}}|\nabla v|^2 dx,
\end{eqnarray}
for any $\alpha\geq 1$ and $\rho$ such that $3\alpha\rho\leq R_2$.
From (\ref{ineqA20}) and using (\ref{ineqA17}) to estimate $N''_w$ in terms of $v$ we get
\begin{equation}
\label{ineqA21}
\int_{B_{2\rho}}|\nabla v|^2 dx\leq C_{10}\left(\frac{\int_{B_{R_1}}|v|^2}{R_1^2\int_{B_{\frac{R_1\lambda_0}{6}}}|\nabla v|^2}\right)^{C_9}  \int_{B_{\rho}}|\nabla v|^2 dx
\end{equation}
for any $\rho\leq \frac{R_1\lambda_0}{6}$ and where $C_9, C_{10}$  depend on $\lambda_0,\lambda_1$ only.\\
Applying (\ref{ineqA21}) to $v-\frac{1}{|B_{R_1}|}\int_{B_{R_1}}v$ and using Poincar\'e inequality we have
\begin{equation}
\label{ineqA22}
\int_{B_{2\rho}}|\nabla v|^2 dx\leq C_{11}\left(\frac{\int_{B_{R_1}}|\nabla v|^2}{R_1^2\int_{B_{\frac{R_1\lambda_0}{6}}}|\nabla v|^2}\right)^{C_9}  \int_{B_{\rho}}|\nabla v|^2 dx
\end{equation}
for $\rho\leq \frac{R_1\lambda_0}{6}$.
Finally we want to prove (\ref{doubeq}). Let $\rho\in [\frac{R_1\lambda_0}{6},\frac{R_0}{2}]$, we have trivially
\begin{eqnarray*}
\int_{B_{2\rho}}|\nabla v|^2 dx&\leq&\int_{B_{R_0}}|\nabla v|^2 dx=\left(\frac{\int_{B_{R_0}}|\nabla v|^2 dx}{\int_{B_{\frac{R_1\lambda_0}{6}}}|\nabla v|^2 dx}\right)\int_{B_{\frac{R_1\lambda_0}{6}}}|\nabla v|^2 dx\\
&\leq &\left(\frac{\int_{B_{R_0}}|\nabla v|^2 dx}{\int_{B_{\frac{R_1\lambda_0}{6}}}|\nabla v|^2 dx}\right)\int_{B_{\rho}}\nabla v|^2 dx.
\end{eqnarray*}
From last inequality and from (\ref{ineqA22})  we immediately get, for $\rho\in [0, R_0/2]$,
\begin{equation}
\label{ineqA23}
\int_{B_{2\rho}}|\nabla v|^2 dx\leq C_{11}\left(\frac{\int_{B_{R_0}}|\nabla v|^2 dx}{\int_{B_{\frac{R_1\lambda_0}{6}}}|\nabla v|^2 dx}\right)\int_{B_{\rho}}|\nabla v|^2 dx.
\end{equation}
Now, we apply the three spheres inequality
\begin{equation}
\label{ineqA24}
\int_{B_{R_0/2}}|\nabla v|^2 dx\leq C_{12}\left(\int_{B_{R_0}}|\nabla v|^2 dx\right)^{\theta}\left(\int_{B_{\frac{R_1\lambda_0}{6}}}|\nabla v|^2 dx\right)^{1-\theta},
\end{equation}
where $ C_{12}$ and $\theta\in (0,1)$ depend on $\lambda_0,\lambda_1$ only. From  (\ref{ineqA24}) we have trivially
\[
\frac{\int_{B_{R_0}}|\nabla v|^2 dx}{\int_{B_{\frac{R_1\lambda_0}{6}}}|\nabla v|^2 dx}\leq C_{12}\left(\frac{\int_{B_{R_0}}|\nabla v|^2 dx}{\int_{B_{\frac{R_0}{2}}}|\nabla v|^2 dx}\right)^{\frac{1}{1-\theta}}.
\]
From last inequality and (\ref{ineqA23}) we finally get (\ref{doubeq}).
\qed
\section*{Acknowledgements}
\noindent
We would like to thank Paolo Bisegna for having proposed the problem to one of the authors pointing out the importance of this problem for  applications. We want also to thank Micol Amar, Daniele Andreucci and Roberto Gianni for stimulating discussions on the problem. This work  was partially supported by Miur grant PRIN 20089PWTPS003.


\begin{thebibliography}{BFV}
\bibitem[AE]{AE} V.~Adolfsson, L.~Escauriaza, $C^{1,\alpha}$ domains and unique continuation at the boundary, \textit{Comm. Pure Appl. Math.}, 50 (1997), 935--969.

\bibitem[ABRV]{ABRV} G.~Alessandrini, E.~Beretta, E.~Rosset, S.~Vessella, Optimal stability for inverse elliptic boundary value problems with unknown boundaries, \textit{Ann. Sc. Norm. Sup. Pisa, Cl. Sci.}, (4) 29 (2000), 755--806.

\bibitem[AMR]{AMR} G.~Alessandrini,  A.~Morassi, E.~Rosset, Size estimates,  in: G. Alessandrini, G. Uhlmann eds., Inverse Problems: Theory and Applications, \textit{Contemporary Mathematics}, 333 (2003), 1--34.

\bibitem[AR]{AR} G.~Alessandrini, E.~Rosset, The inverse conductivity problem with one measurement: bounds on the size of the unknown object, \textit{SIAM J. Appl. Math.}, 58, no. 4 (1998), 1060--1071.

 \bibitem[ARS]{ARS} G.~Alessandrini, E.~Rosset,  J.K.~Seo, Optimal size estimates for the inverse conductivity problem with one measurement, \textit{Proc. Amer. Math. Soc.}, 128 (2000), 53--64.


 \bibitem[ARRoV]{ARRoV} G.~Alessandrini, L.~Rondi, E.~Rosset, S.~Vessella, The stability for the Cauchy problem for elliptic equations , \textit{Inverse Problems}, 25 (2009), 1--47.


\bibitem[Bo]{Bo} L.~Borcea, Electrical impedance tomography, \textit{Inverse Problems}, 18, no. 6 (2002), R99--R136.


\bibitem[CIN]{CIN} M.~Cheney, D.~Isaacson, J. C.~Newell, Electrical Impedence Tomography, \textit{Siam Review} 41, No. 1, (1999), 85--101.


\bibitem[CG]{CG} A.V.~Cherkaev,  L. V. ~Gibiansky, Variational principles for complex conductivity, viscoelasticiy, and similar problems in media with complex moduli, \textit{J. Math. Phys.}, 35 (1994), 127--145.



\bibitem[CV]{CV} Y.~Capdeboscq,  M.~Vogelius, Optimal asymptotic estimates for the volume of internal
inhomogeneities in terms of multiple boundary measurements, \textit{Math. Modelling Num. Anal.}, 37 (2003), 227--240.

\bibitem[CoFe]{CoFe} R.~Coifman, C.~Fefferman, Weighted norm inequalities for maximal functions and singular integrals, \textit{Studia Math.}, 51 (1974), 241--250.

\bibitem[GC-RDeF]{GC-RDeF} J.~Garc\'{\i}a-Cuerva, J.L.~Rubio De Francia, "Weighted norm inequalities and related topics", Amsterdam, North Holland, 1985.

\bibitem[GL]{GL} N.~Garofalo, F.H.~Lin, Monotonicity properties of variational integrals, $A_p$ weights and unique continuation, \textit{Indiana Univ. Math. J.}, 35 (1986), 245--268.

\bibitem[GT]{GT} D.~Gilbarg, N.S.~Trudinger, "Elliptic Partial Differential Equations of Second Order", Grundleheren der Mathematischen Wissenshaften vol.224, 2nd edition, Berlin, Springer, 1983.


\bibitem[G]{G} H.~Griffiths, Tissue Spectroscopy with Electrical Impedance
Tomography: Computer Simulations, \textit{IEEE Transactions on Biomedical Engineering} 42, No.9 (1995) 948--953

\bibitem[LNW]{LNW}C.L.~Lin, G.~Nakamura, J.N.~Wang, Quantitative uniqueness for second order elliptic operators
with strongly singular coefficients, \textit{Rev. Mat. Iberoamericana}, 27, n. 2 (2011), 475--491.

\bibitem[KKM]{KKM} H.~Kang, E. ~Kim, G.~Milton, Sharp bounds on the volume fractions of two materials in a
two-dimensional body from electrical boundary measurements: the translation method, \textit{ arXiv preprint},  (2011).

\bibitem[KSS]{KSS} H.~Kang, J.K.~Seo, D.~Sheen, The inverse conductivity problem with one measurement:
stability and estimation of size, \textit{SIAM J. Math. Anal.}, 28 (1997), 1389--1405.

\bibitem[Ku]{Ku} I.~Kukavica, Quantitative uniqueness for second order elliptic operators, \textit{Duke Math. J.}, 91 (19..), 225-240.

\bibitem[Ve]{Ve} S.~Vessella, Quantitative continuation fron a measurable set of solutions of elliptic equations, \textit{Proc R. Soc. Edinb.}, 130 (A), 909--923.
\end{thebibliography}
\end{document}